\documentclass[11pt,reqno]{amsart}

\usepackage{amsmath,amssymb,mathrsfs}
\usepackage{graphicx,cite,cases,arydshln}
\usepackage[]{changes}

\numberwithin{equation}{section}

\begin{document}

\title[a fast direct imaging method]{A fast direct imaging method for the
inverse obstacle scattering problem with nonlinear point scatterers}

\author{Jun Lai}
\address{School of Mathematical Sciences, Zhejiang
University, Hangzhou, Zhejiang 310027, China}
\email{laijun6@zju.edu.cn}

\author{Ming Li}
\address{College of Data Science, Taiyuan University of Technology, Taiyuan,
Shangxi 030024, China}
\email{liming01@tyut.edu.cn}

\author{Peijun Li}
\address{Department of Mathematics, Purdue University, West Lafayette, IN 47907,
USA.}
\email{lipeijun@math.purdue.edu}

\author{Wei Li}
\address{Institute for Mathematics and its Applications, Minneapolis, MN 55455,
USA.}
\email{lixx1774@umn.edu}

\thanks{The work of J. Lai was supported in part by the Funds for Creative
Research Groups of NSFC (No. 11621101) and the Major Research Plan of NSFC (No.
91630309). The research of M. Li was supported partially by the National Youth
Science Foundation of China (Grant no. 11401423). The research of P. Li was
supported in part by the NSF grant DMS-1151308.}

\subjclass[2010]{78A46, 78M15, 65N21}

\keywords{Foldy--Lax formulation, point scatterers, inverse obstacle scattering
problem, the Helmholtz equation, boundary integral equation, nonlinear optics}

\begin{abstract}
Consider the scattering of a time-harmonic plane wave by heterogeneous media
consisting of linear or nonlinear point scatterers and extended obstacles. A
generalized Foldy--Lax formulation is developed to take fully into account of
the multiple scattering by the complex media. A new imaging function is proposed
and an FFT-based direct imaging method is developed for the inverse
obstacle scattering problem, which is to reconstruct the shape of the
extended obstacles. The novel idea is to utilize the nonlinear point scatterers
to excite high harmonic generation so that enhanced imaging resolution can be
achieved. Numerical experiments are presented to demonstrate the effectiveness
of the proposed method.
\end{abstract}

\maketitle

\section{Introduction}

In scattering theory, one of the basic problems is the scattering of
a time-harmonic plane wave by an impenetrable medium, which is called the
obstacle scattering problem \cite{ck-83}. Given the incident wave, the direct
obstacle scattering problem is to determine the scattered wave for the known
obstacle; while the inverse obstacle scattering problem is to reconstruct the
shape of the obstacle from either near-field or far-field measurement of the
scattered wave. The obstacle scattering problem has played a fundamental role in
diverse scientific areas such as geophysical exploration, radar and sonar,
medical imaging, and nondestructive testing. 

The inverse obstacle scattering problem is challenging due to nonlinearity and
ill-posedness. It has been extensively studied by many researchers. Various
computational methods have been developed to overcome the issues and solve the
inverse problem. Broadly speaking, these methods can be classified into two
types: optimization based iterative methods \cite{bhl-jcp07} and imaging based
direct methods \cite{agkls-sjis12, agm-sjis14, akkll-sinum11, cc-06, c-ip01,
cln-ip05, ck-ip96, d, ep-sjsc06, gmd-jasa04, hsz-jcp04, hsz-ip06, hsz-ip07,
hhsz-jasa09, i-ip98, kpf-jasa03, kg-08, p-jcam00}. The former are known as
quantitative methods and aim at recovering the functions
which parameterize the obstacles. The latter are usually called qualitative
methods and attempt to design imaging functions which highlight the obstacles.
According to the Rayleigh criterion, there is a resolution limit, roughly one
half the wavelength, on the accuracy of the reconstruction for a given incident
wave \cite{bpz-jasa02, bgpt-ip11, z-siap04}. To improve the resolution, one
approach is simply to use an incident wave with shorter wavelength or higher
frequency as an illumination. A topical review can be found in \cite{bllt-ip15}
on computational approaches and mathematical analysis for solving inverse
scattering problems by multi-frequencies. 

In this paper, we consider the inverse obstacle scattering problem with
an emphasis on resolution enhancement. Motivated by nonlinear optics
\cite{blms-17, b-08, ls-pr15}, we utilize nonlinear point scatterers to excite
high harmonic generation so that enhanced resolution can be achieved to
reconstruct the obstacles. To realize this idea, we have to consider the
scattering problem of a time-harmonic plane incident wave by a heterogeneous
medium, which consists of small scale point scatterers and wavelength comparable
extended obstacles. The main purpose of this work is threefold: 

\begin{enumerate}
 
 \item develop a generalized Foldy--Lax formulation and design an efficient
solver for the scattering problem involving point scatterers and
extended obstacles; 
 
 \item propose a new imaging function and develop an FFT-based method to
efficiently evaluate the imaging function on large sampling points; 

 \item explore the features of high harmonic generation for nonlinear optics and
apply them to the area of inverse scattering to achieve enhanced imaging
resolution. 

\end{enumerate}

\begin{figure}
\centering
\includegraphics[width=0.35\textwidth]{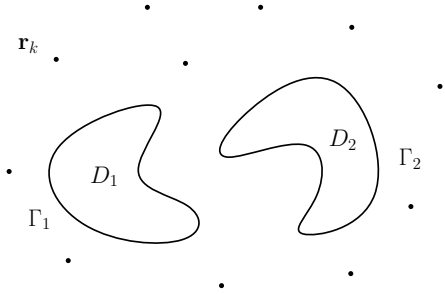}
\caption{Schematic of the problem geometry.}
\label{pg}
\end{figure}

The Foldy--Lax formulation was developed in \cite{f-pr45, l-rmp51} to
describe the scattering of an incoming wave by a group of linear
point scatterers. The scattered field can be computed by solving a
self-contained system of liner equations. Using a unified approach, we first
extend the Foldy--Lax formulation to handle a group of linear, quadratically
nonlinear, or cubically nonlinear point scatterers. It is known that the
Foldy--Lax formulation is only appropriate for media whose sizes are much
smaller than the wavelength \cite{chs-m3as14, cs-ip12, cs-mms14,
vcl-rmp98, m-06}. It is no longer adequate for the scattering by wavelength
comparable media \cite{hms-mms14, llz-mms14}. The boundary integral equation
method is particularly useful for the scattering by extended obstacles. The
scattering problem becomes sophisticated when both the point scatterers and
extended obstacles are present, as seen in Figure \ref{pg}. The generalized
Foldy--Lax formulation has been studied in \cite{bhlz-mc14, hl-mms10, hlz-jcp13,
hsz-jcp10} to take into account the multiple scattering between linear point
scatterers and extended obstacles. Here we develop a generalized Foldy--Lax
formulation to fully take into account the multiple scattering between the
nonlinear point scatterers and extended obstacles. The generalized formulation
couples the Foldy--Lax and boundary integral equation formulations and is
self-contained linear or nonlinear system of equations. The coupled system needs
to be solved numerically in order to obtain the scattered field and its
far-field pattern. For linear point scatterers, we apply LU factorization based
direct solver to solve the coupled linear system of equations; for nonlinear
point scatterers, we propose an efficient nonlinear solver which combines the
Schur complement technique and trust region Newton type method. 

The imaging based methods do not require solving any direct problem, which
make them very attractive to solve the inverse scattering problem. But, they
could be still very time-consuming when evaluating the imaging functions on
large sampling points. The reason is that the evaluation usually involves the
dense matrix-vector multiplication at each sampling point. To overcome this
issue, we construct two illumination vectors and propose a new imaging function.
The function, under an appropriate transformation, can be taken as the Fourier
transform of the response matrix from the frequency space into the physical
space. However, the frequency sampling points are not uniform, which implies the
standard FFT can not be applied. We propose an acceleration technique based on
the non-uniform fast Fourier transform (NUFFT)\cite{gl-sr04}, which highly
reduces the computational complexity and accelerates the evaluation. 

Using the quadratically or cubically nonlinear point scatterers, the second or
third harmonic generation is excited to image the extended obstacles.
Essentially, this approach is equivalent to using double or
triple frequency wave to illuminate the obstacles. As a consequence, enhanced
resolution can be achieved. However, we find out that the location of the
nonlinear point scatterers is crucial to the imaging. Interestingly, they should
be aligned up in the same direction as the propagation direction of the plane
incident wave so that the correct reconstruction results can appear with the
desirable imaging resolution. Numerical experiments are shown to demonstrate the
effectiveness of the proposed method. To make the paper self-contained, we
briefly introduce in the appendix the nonlinear wave equations and nonlinear
point scatterer models. 

The paper is organized as follows. In section 2, the Foldy--Lax formulation is
introduced for the scattering by a group of point scatterers. In section 3, the
boundary integral formulation is briefly reviewed for the scattering by
extended scatterers. Section 4 introduces the Generalized Foldy--Lax
formulation for the scattering by mixed scatterers. Section 5 is devoted to the
inverse scattering problem, where the fast direct imaging method is developed.
Numerical experiments are given in section 6. The paper is concluded in section
7.

\section{Foldy--Lax formulation}

In this section, we introduce the Foldy--Lax formulation for the scattering of
a plane incident wave by a group of linear, quadratically nonlinear, and
cubically nonlinear point scatterers. 

\subsection{Linear point scatterers}

Consider a collection of $m$ separated linear and isotropic point scatterers
located at ${\bf r}_k\in\mathbb R^2, k=1, \dots, m$. Let $\phi_{\rm inc}$ be the
plane incident wave given explicitly by 
\begin{equation}\label{if}
 \phi_{\rm inc}({\bf r})=e^{{\rm i}\kappa {\bf r}\cdot{\bf d}}, \quad {\bf
r}\in\mathbb R^2, 
\end{equation}
where $\kappa=\omega/c$ is the wavenumber, $\omega>0$ is the angular frequency
and $c$ is the speed of wave propagation, and ${\bf d}\in\mathbb S^1$ is the
unit propagation direction. It is easy to verify the incident field
satisfies the Helmholtz equation:
\begin{equation}\label{he-if}
 \Delta\phi_{\rm inc}+\kappa^2\phi_{\rm inc}=0\quad\text{in}~\mathbb R^2. 
\end{equation}

As is shown in \eqref{b-lp}, the total field $\phi$ satisfies
\begin{equation}\label{flhe-tf}
\Delta\phi({\bf r})+\kappa^2 \phi({\bf
r})=-\sum_{k=1}^m\sigma_k\phi_k\delta({\bf r}-{\bf r}_k),\quad{\bf r}\in\mathbb
R^2,
\end{equation}
where $\sigma_k>0$ is the scattering coefficient for the $k$-th point
scatterer, $\phi_k$ is the external field acting on the $k$-th point
scatterer, and $\delta$ is the Dirac delta function. We comment that $\phi_k$
is also the total field at the location ${\bf r}_k$, i.e., $\phi_k=\phi({\bf
r}_k)$. 

The total field $\phi$ consists of the incident field $\phi_{\rm inc}$ and the
scattered field $\psi$: 
\[
 \phi=\phi_{\rm inc}+\psi.
\]
Subtracting \eqref{he-if} from \eqref{flhe-tf} yields 
\begin{equation}\label{flhe-sf}
\Delta\psi({\bf r})+\kappa^2 \psi({\bf
r})=-\sum_{k=1}^m\sigma_k\phi_k\delta({\bf r}-{\bf r}_k),\quad{\bf r}\in\mathbb
R^2.
\end{equation}
The scattered field is required to satisfy the Sommerfeld radiation condition:
\begin{equation}\label{fl-rc}
 \lim_{r\to\infty}r^{1/2}(\partial_r\psi-{\rm i}\kappa\psi)=0,\quad r=|{\bf
r}|. 
\end{equation}

It follows from \eqref{flhe-sf} and \eqref{fl-rc} that the scattered field can
be written as 
\begin{equation}\label{fl-lsf}
\psi({\bf r})=\sum_{k=1}^m \sigma_k\phi_k G_\kappa({\bf r}, {\bf r}_k), 
\end{equation}
where 
\[
 G_\kappa({\bf r}, {\bf r}')=\frac{\rm i}{4}H_0^{(1)}(\kappa|{\bf r}-{\bf r}'|)
\]
is the free space Green's function for the two-dimensional
Helmholtz equation. Here $H_0^{(1)}$ is the Hankel function of the first kind
with order zero. It is left to compute the external fields $\phi_k$ in order to
compute the scattered field \eqref{fl-lsf}. 

Adding the incident field on both sides of \eqref{fl-lsf} gives
\begin{equation}\label{fl-ltf}
 \phi({\bf r})=\phi_{\rm inc}({\bf r})+\sum_{k=1}^m \sigma_k\phi_k G_\kappa({\bf
r}, {\bf r}_k). 
\end{equation}
 Evaluating \eqref{fl-ltf} at ${\bf r}_i$ and excluding the self-interaction to
avoid the singularity of the Green function, we obtain a linear system of
algebraic equations for $\phi_k$:
\begin{equation}\label{fl-l}
 \phi_i - \sum_{\stackrel{\scriptstyle k=1}{k\ne
i}}^{m} \sigma_k \phi_k G_\kappa({\bf r}_i, {\bf r}_k)=\phi_{\rm inc}({\bf
r}_i),
\end{equation}
which is known as the Foldy--Lax formulation. 

\subsection{Quadratically nonlinear point scatterers}

We assume that the nonlinear point scatterers respond to the incoming wave
quadratically and the nonlinearity is weak. Let the plane incident wave
$\phi_{\rm inc}$, given in \eqref{if}, be of a single frequency $\omega$. The
point scatterers generate fields at frequencies $\omega_1=\omega$ and
$\omega_2=2\omega$ due to the quadratic nonlinearity, which is known as the
second harmonic generation. Let $\kappa_j=\omega_j/c$ and $\phi^{(j)}$ be the
field at the frequency $\omega_j, j=1, 2$.

As is known in \eqref{b-qp}, the fields $\phi^{(j)}$ satisfy the coupled
Helmholtz equations in $\mathbb R^2$: 
\begin{subequations}\label{flhe-qtf}
 \begin{align}
  \Delta \phi^{(1)}({\bf r})+\kappa_1^2 \phi^{(1)}({\bf r})&=-\sum_{k=1}^m
\left(\sigma_{k,1}^{(1)}\phi_k^{(1)}+\sigma_{k,1}^{(2)}\bar{\phi}_k^{(1)}\phi_k^{(2)}
\right)\delta( { \bf r}-{\bf r}_k),\\
 \Delta \phi^{(2)}({\bf r})+\kappa_2^2 \phi^{(2)}({\bf r})&=-\sum_{k=1}^m
\left(\sigma_{k,2}^{(1)}\phi_k^{(2)}+\sigma_{k,2}^{(2)}\bigl(\phi_k^{(1)}\bigr)^2
\right)\delta({ \bf r}-{\bf r}_k),
 \end{align}
\end{subequations}
where the bar denotes the complex conjugate, $\sigma_{k,1}^{(1)}$ and
$\sigma_{k,2}^{(1)}$ are the linear scattering coefficients for the $k$-th point
scatterer, $\sigma_{k,1}^{(2)}$ and $\sigma_{k,2}^{(2)}$ are the quadratically
nonlinear scattering coefficients for the $k$-th point scatterer, and
$\phi_k^{(j)}$ is the external field acting on the $k$-th point scatterer at
frequency $\omega_j$. The fields $\phi^{(j)}$ satisfy the following
relationship: 
\[
 \phi^{(1)}=\phi_{\rm inc}+\psi^{(1)},\quad \phi^{(2)}=\psi^{(2)},
\]
where $\psi^{(j)}$ is the scattered field corresponding to the wavenumber
$\kappa_j$. It can be verified that the scattered fields satisfy  
\begin{subequations}\label{flhe-qsf}
 \begin{align}
  \Delta \psi^{(1)}({\bf r})+\kappa_1^2 \psi^{(1)}({\bf r})&=-\sum_{k=1}^m
\left(\sigma_{k,1}^{(1)}\phi_k^{(1)}+\sigma_{k,1}^{(2)}\bar{\phi}_k^{(1)}\phi_k^{(2)}
\right)\delta( { \bf r}-{\bf r}_k),\\
 \Delta \psi^{(2)}({\bf r})+\kappa_2^2 \psi^{(2)}({\bf r})&=-\sum_{k=1}^m
\left(\sigma_{k,2}^{(1)}\phi_k^{(2)}+\sigma_{k,2}^{(2)}\bigl(\phi_k^{(1)}\bigr)^2
\right)\delta({ \bf r}-{\bf r}_k).
 \end{align}
\end{subequations}
In addition, they are required to satisfy the Sommerfeld radiation condition
\[
 \lim_{r\to\infty}r^{1/2}\left(\partial_r\psi^{(j)}-{\rm
i}\kappa_j\psi^{(j)}\right)=0,\quad r=|{\bf r}|.
\]

It follows from \eqref{flhe-qsf} that the scattered fields satisfy 
\begin{subequations}\label{fl-qsf}
 \begin{align}
\label{fl-qsf1}\psi^{(1)}({\bf r})&=
\sum_{k=1}^m \left(
\sigma_{k,1}^{(1)}\phi_k^{(1)}+\sigma_{k,1}^{(2)}\bar{\phi}_k^{(1)}\phi_k^{(2)}
\right) G_{\kappa_1}({\bf r}, {\bf r}_k),\\
\label{fl-qsf2}\psi^{(2)}({\bf r})&=\sum_{k=1}^m
\left(\sigma_{k,2}^{(1)}\phi_k^{(2)}+\sigma_{k,2}^{(2)}\bigl(\phi_k^{(1)}\bigr)^2\right)
G_{\kappa_2}({\bf r}, {\bf r}_k),
 \end{align}
\end{subequations}
where $G_{\kappa_j}$ is the free space Green's function for the two-dimensional
Helmholtz equation at the wavenumber $\kappa_j$. Adding the incident field on
both sides of \eqref{fl-qsf1} and noting $\phi^{(2)}=\psi^{(2)}$ for
\eqref{fl-qsf2}, we obtain
\begin{subequations}\label{fl-qtf}
 \begin{align}
  \phi^{(1)}({\bf r})&=\phi_{\rm inc}({\bf r})+\sum_{k=1}^m
\left(\sigma_{k,1}^{(1)}\phi_k^{(1)}+\sigma_{k,1}^{(2)}\bar{\phi}_k^{(1)}\phi_k^{(2)}
\right) G_{\kappa_1}({\bf r}, {\bf r}_k),\\
\phi^{(2)}({\bf r})&=\sum_{k=1}^m
\left(\sigma_{k,2}^{(1)}\phi_k^{(2)}+\sigma_{k,2}^{(2)}\bigl(\phi_k^{(1)}\bigr)^2\right)
G_{\kappa_2}({\bf r}, {\bf r}_k). 
 \end{align}
\end{subequations}
Similarly, evaluating \eqref{fl-qtf} at ${\bf r}_i$ and excluding the
self-interaction yield a nonlinear system of equations for $\phi_k^{(j)}$:
\begin{subequations}\label{fl-q}
 \begin{align}
   \phi^{(1)}_i&-\sum_{\stackrel{\scriptstyle k=1}{k\ne i}}^{m}
\left(\sigma_{k,1}^{(1)}\phi_k^{(1)}+\sigma_{k,1}^{(2)}\bar{\phi}_k^{(1)}\phi_k^{(2)}
\right) G_{\kappa_1}({\bf r}_i, {\bf r}_k)=\phi_{\rm inc}({\bf r}_i),\\
\phi^{(2)}_i&-\sum_{\stackrel{\scriptstyle k=1}{k\ne
i}}^{m} \left(
\sigma_{k,2}^{(1)}\phi_k^{(2)}+\sigma_{k,2}^{(2)}\bigl(\phi_k^{(1)}\bigr)^2
\right) G_{\kappa_2}({\bf r}_i, {\bf r}_k)=0,
 \end{align}
\end{subequations}
which is the Foldy--Lax formulation for point scatterers with quadratic
nonlinearity.

\subsection{Cubically nonlinear point scatterers}

Taking the plane incident wave $\phi_{\rm inc}$ with frequency $\omega$, we
consider the scattering by point scatterers with weak cubic nonlinearity. The
interaction gives rise to fields with frequencies $\omega_1=\omega$ and
$\omega_3=3\omega$, which is called the third harmonic generation.
Let $\kappa_j=\omega_j/c$ and denote the field at frequency $\omega_j$ by
$\phi^{(j)}, j=1, 3$. 

By \eqref{b-cp}, the fields $\phi^{(j)}$ satisfy the following coupled Helmholtz
equations in $\mathbb R^2$: 
\begin{subequations}\label{flhe-ctf}
 \begin{align}
  \Delta \phi^{(1)}({\bf r})+\kappa_1^2 \phi^{(1)}({\bf r})&=-\sum_{k=1}^m
\left( \sigma_{k,1}^{(1)}\phi^{(1)}_k+\sigma_{k,1}^{(3)}|\phi_k^{(1)}|^2\phi_k^{
(1)}+\sigma_{k,2}^{(3)}\bigl(\bar{\phi}_k^{(1)}\bigr)^2\phi_k^{(3)}
\right)\delta( { \bf r}-{\bf r}_k),\\
 \Delta \phi^{(3)}({\bf r})+\kappa_2^2 \phi^{(3)}({\bf
r})&=-\sum_{k=1}^m\left( \sigma_{k,2}^{(1)}\phi_k^{(3)}+\sigma_{k,3}^{(3)}
\bigl(\phi_k^{(1)}\bigr)^3 \right) \delta({\bf r}-{\bf r}_k),
\end{align}
\end{subequations}
where $\sigma_{k,1}^{(1)}$ and $\sigma_{k,2}^{(1)}$ are the linear
scattering coefficients for the $k$-th point scatterer, $\sigma_{k,1}^{(3)}$,
$\sigma_{k,2}^{(3)}$ and $\sigma_{k,3}^{(3)}$ are the cubically 
nonlinear scattering coefficients for the $k$-th point scatterer, and
$\phi_k^{(j)}$ is the external field acting on the $k$-th point scatterer at
frequency $\omega_j$. The fields $\phi^{(j)}$ satisfy the following
relationship: 
\[
 \phi^{(1)}=\phi_{\rm inc}+\psi^{(1)},\quad \phi^{(3)}=\psi^{(3)},
\]
where $\psi^{(j)}$ is the scattered field corresponding to the wavenumber
$\kappa_j$. It can be verified that the scattered fields satisfy  
\begin{subequations}\label{flhe-csf}
 \begin{align}
  \Delta \psi^{(1)}({\bf r})+\kappa_1^2 \psi^{(1)}({\bf r})&=-\sum_{k=1}^m
\left( \sigma_{k,1}^{(1)}\phi^{(1)}_k+\sigma_{k,1}^{(3)}|\phi_k^{(1)}|^2\phi_k^{
(1)}+\sigma_{k,2}^{(3)}\bigl(\bar{\phi}_k^{(1)}\bigr)^2\phi_k^{(3)}
\right)\delta({\bf r}-{\bf r}_k),\\
 \Delta \psi^{(3)}({\bf r})+\kappa_2^2 \psi^{(3)}({\bf
r})&=-\sum_{k=1}^m\left( \sigma_{k,2}^{(1)}\phi_k^{(3)}+\sigma_{k,3}^{(3)}
\bigl(\phi_k^{(1)}\bigr)^3\right)\delta({\bf r}-{\bf r}_k).
 \end{align}
\end{subequations}
In addition, they are required to satisfy the Sommerfeld radiation condition
\[
 \lim_{r\to\infty}r^{1/2}\left(\partial_r\psi^{(j)}-{\rm
i}\kappa_j\psi^{(j)}\right)=0,\quad r=|{\bf r}|.
\]

It is easy to verify from \eqref{flhe-csf} that the scattered fields satisfy 
\begin{subequations}\label{fl-csf}
 \begin{align}
\label{fl-csf1}\psi^{(1)}({\bf r})&=\sum_{k=1}^m
\left( \sigma_{k,1}^{(1)}\phi^{(1)}_k+\sigma_{k,1}^{(3)}|\phi_k^{(1)}|^2\phi_k^{
(1)}+\sigma_{k,2}^{(3)}\bigl(\bar{\phi}_k^{(1)}\bigr)^2\phi_k^{(3)}
\right) G_{\kappa_1}({\bf r}, {\bf r}_k),\\
\label{fl-csf2}\psi^{(3)}({\bf r})&=\sum_{k=1}^m\left(
\sigma_{k,2}^{(1)}\phi_k^{(3)}+\sigma_{k,3}^{(3)} \bigl(\phi_k^{(1)}\bigr)^3
\right) G_{\kappa_3}({\bf r}, {\bf r}_k),
 \end{align}
\end{subequations}
where $G_{\kappa_j}$ is the free space Green's function for the two-dimensional
Helmholtz equation at the wavenumber $\kappa_j$. Adding the incident field on
both sides of \eqref{fl-csf1} and noting $\phi^{(2)}=\psi^{(2)}$ for
\eqref{fl-csf2}, we obtain
\begin{subequations}\label{fl-ctf}
 \begin{align}
  \phi^{(1)}({\bf r})&=\phi_{\rm inc}({\bf r})+\sum_{k=1}^m
\left(\sigma_{k,1}^{(1)}\phi^{(1)}_k+\sigma_{k,1}^{(3)}|\phi_k^{(1)}|^2\phi_k^{
(1)}+\sigma_{k,2}^{(3)}\bigl(\bar{\phi}_k^{(1)}\bigr)^2\phi_k^{(3)}
\right)G_{\kappa_1}({\bf r}, {\bf r}_k),\\
\phi^{(3)}({\bf r})&=\sum_{k=1}^m\left(
\sigma_{k,2}^{(1)}\phi_k^{(3)}+\sigma_{k,3}^{(3)} \bigl(\phi_k^{(1)}\bigr)^3
\right) G_{\kappa_2}({\bf r}, {\bf r}_k). 
 \end{align}
\end{subequations}
Evaluating \eqref{fl-ctf} at ${\bf r}_i$ and excluding the
self-interaction, we get a nonlinear system of equations for $\phi_k^{(j)}$:
\begin{subequations}\label{gl-c}
 \begin{align}
   \phi^{(1)}_i&-\sum_{\stackrel{\scriptstyle k=1}{k\ne i}}^m
\left( \sigma_{k,1}^{(1)}\phi^{(1)}_k+\sigma_{k,1}^{(3)}|\phi_k^{(1)}|^2\phi_k^{
(1)}+\sigma_{k,2}^{(3)}\bigl(\bar{\phi}_k^{(1)}\bigr)^2\phi_k^{(3)}
\right)G_{\kappa_1}({\bf r}_i, {\bf r}_k)=\phi_{\rm inc}({\bf r}_i),\\
\phi^{(3)}_i&-\sum_{\stackrel{\scriptstyle k=1}{k\ne i}}^m\left(
\sigma_{k,2}^{(1)}\phi_k^{(3)}+\sigma_{k,3}^{(3)}\bigl(\phi_k^{(1)}\bigr)^3
\right) G_{\kappa_3}({\bf r}_i, {\bf r}_k)=0,
 \end{align}
\end{subequations}
which is the Foldy--Lax formulation for point scatterers with cubic
nonlinearity.

\section{Boundary integral formulation}

In this section, we briefly introduce the boundary integral equation method for
solving the scattering problem with extended scatterers. The detailed
information can be found in \cite{ck-83}. 

Consider the scattering of a plane incident wave by a sound-soft extended
scatterer, which is described by the domain $D$ with a boundary $\Gamma$
consisting of a finite number of disjoint, closed, bounded surfaces belonging to
the class $C^2$. The exterior $\mathbb R^2\setminus\bar{D}$ is assumed to the
connected, whereas $D$ itself is allowed to have more than one component. The
unit normal vector $\nu$ on $\Gamma$ is assumed to be directed into the exterior
of $D$.

The total field satisfies the Helmholtz equation:
\begin{equation}\label{bihe-tf}
 \Delta\phi+\kappa^2\phi=0\quad\text{in}~ \mathbb R^2\setminus\bar{D}.
\end{equation}
The sound-soft boundary implies that 
\begin{equation}\label{bihe-bc}
 \phi=0\quad\text{on} ~ \Gamma.
\end{equation}
The total field $\phi$ consists of the incident field $\phi_{\rm inc}$ and the
scattered field $\psi$:
\begin{equation}\label{bihe-fd}
 \phi=\phi_{\rm inc}+\psi.
\end{equation}
It follows from \eqref{he-if}, \eqref{bihe-tf}, and \eqref{bihe-fd} that the
scattered field $\psi$ also satisfies the Helmholtz equation:
\begin{equation}\label{bihe-sf}
 \Delta\psi+\kappa^2\psi=0\quad\text{in}~ \mathbb R^2\setminus\bar{D}.
\end{equation}
The following Sommerfeld radiation condition is imposed to ensure the
well-posedness of the scattering problem:
\begin{equation}\label{bi-rc}
 \lim_{r\to\infty}r^{1/2}(\partial_r\psi-{\rm i}\kappa\psi)=0,\quad r=|{\bf
r}|. 
\end{equation}

It is shown in \cite{ck-83} by using the second Green theorem that 
\[
 \int_\Gamma \left(\phi_{\rm inc}({\bf r}')\partial_{\nu'} G_\kappa({\bf r},
{\bf r}')-\partial_{\nu'}\phi_{\rm inc}({\bf r}')G_\kappa({\bf r}, {\bf
r}')\right) {\rm d}s({\bf r}')=\begin{cases}
                                -\phi_{\rm inc}({\bf r}),&\quad{\bf r}\in D,\\
                                0,&\quad{\bf r}\in\mathbb R^2\setminus\bar{D}.
                               \end{cases}
\]
and
\[
 \int_\Gamma \left(\psi({\bf r}')\partial_{\nu'} G_\kappa({\bf r},
{\bf r}')-\partial_{\nu'}\psi({\bf r}')G_\kappa({\bf r}, {\bf
r}')\right) {\rm d}s({\bf r}')=\begin{cases}
                                0,&\quad{\bf r}\in D,\\
                                \psi({\bf r}),&\quad{\bf r}\in\mathbb
R^2\setminus\bar{D}.
                               \end{cases}
\]
Adding the above two equations and using the sound-soft boundary
condition \eqref{bihe-bc}, we get 
\begin{equation}\label{bi-if}
 \phi_{\rm inc}({\bf r})=\int_\Gamma G_\kappa({\bf r}, {\bf
r}')\partial_{\nu'}\phi({\bf r}'){\rm d}s({\bf r}'),\quad {\bf r}\in D.
\end{equation}
and
\begin{equation}\label{bi-sf}
 \psi({\bf r})=-\int_\Gamma G_\kappa({\bf r}, {\bf
r}')\partial_{\nu'}\phi({\bf r}'){\rm d}s({\bf r}'),\quad{\bf r}\in\mathbb
R^2\setminus\bar{D}.
\end{equation}
To compute the scattered field $\psi$, it is required to compute
$\partial_\nu\phi$ on $\Gamma$.

Taking the normal derivative of \eqref{bi-if} on $\Gamma$ and applying the jump
condition yield
\begin{equation}\label{bi-nif}
 \partial_\nu\phi_{\rm inc}({\bf r})=\int_\Gamma \partial_\nu G_\kappa({\bf r},
{\bf r}')\partial_{\nu'}\phi({\bf r}'){\rm d}s({\bf
r}')+\frac{1}{2}\partial_\nu\phi({\bf r}),\quad {\bf r}\in\Gamma.
\end{equation}
Multiplying \eqref{bi-if} by ${\rm i}\eta$ and subtracting it from
\eqref{bi-nif}, we thus obtain a boundary integral equation for
$\partial_\nu\phi$ on $\Gamma$:
\begin{equation}\label{bi}
 \frac{1}{2}\partial_\nu\phi({\bf r})+\int_\Gamma (\partial_\nu-{\rm
i}\eta) G_\kappa({\bf r}, {\bf r}')\partial_{\nu'}\phi({\bf r}'){\rm
d}s({\bf r}')=(\partial_\nu-{\rm i}\eta)\phi_{\rm inc}({\bf r}),
\end{equation}
where the coupling parameter $\eta>0$ is introduced to ensure the unique
solvability of \eqref{bi}.

\section{Generalized Foldy--Lax formulation}

This section presents the generalized Foldy--Lax formulation for the scattering
problem of mixed scatterers, which consist of both the extended and point
scatterers. We introduce the generalized Foldy--Lax formulation for the linear,
quadratically nonlinear, and cubically nonlinear point scatterers,
respectively. 

\subsection{Linear point scatterers}

Viewing the external field acting on the point scatterers as point sources for
the extended obstacle, we have the following equation for the total field:
\begin{equation}\label{gflhe-tf}
\Delta \phi({\bf r}) + \kappa^2 \phi({\bf r}) = -\sum_{k=1}^m \sigma_k \phi_k
\delta({\bf r}-{\bf r}_k),\quad{\bf r}\in\mathbb{R}^2\setminus\bar{D},
\end{equation}
where $\phi_k$ is the external field acting on the $k$-th point scatterer and
$\delta$ is the Dirac delta function. The obstacle is still assumed to be
sound-soft. The total field vanishes on the boundary, i.e., 
\begin{equation}\label{gflhe-bc}
 \phi=0\quad\text{on}~ \Gamma. 
\end{equation}
Subtracting the incident field \eqref{he-if} from the total field
\eqref{gflhe-tf}, we get the equation for the scattered field:
\begin{equation}\label{gflhe-sf}
\Delta \psi({\bf r}) + \kappa^2 \psi({\bf r}) = -\sum_{k=1}^m \sigma_k \phi_k
\delta({\bf r}-{\bf r}_k),\quad{\bf r}\in\mathbb{R}^2\setminus\bar{D},
\end{equation}
As usual, the scattered field is required to satisfy the Sommerfeld radiation
condition:
\begin{equation}\label{gfl-rc}
 \lim_{r\to\infty}r^{1/2}(\partial_r\psi-{\rm i}\kappa\psi)=0,\quad r=|{\bf
r}|. 
\end{equation}

We can follow the same steps as those in \cite{ck-83} to show that 
\[
 \int_\Gamma \left(\phi_{\rm inc}({\bf r}')\partial_{\nu'} G_\kappa({\bf r},
{\bf r}')-\partial_{\nu'}\phi_{\rm inc}({\bf r}')G_\kappa({\bf r}, {\bf
r}')\right) {\rm d}s({\bf r}')=\begin{cases}
                                -\phi_{\rm inc}({\bf r}),&\quad{\bf r}\in D,\\
                                0,&\quad{\bf r}\in\mathbb R^2\setminus\bar{D}.
                               \end{cases}
\]
and
\[
 \sum_{k=1}^m \sigma_k\phi_k G_\kappa({\bf r}, {\bf r}_k)+\int_\Gamma
\left(\psi({\bf r}')\partial_{\nu'} G_\kappa({\bf r},
{\bf r}')-\partial_{\nu'}\psi({\bf r}')G_\kappa({\bf r}, {\bf
r}')\right) {\rm d}s({\bf r}')=\begin{cases}
                                0,&\quad{\bf r}\in D,\\
                                \psi({\bf r}),&\quad{\bf r}\in\mathbb
R^2\setminus\bar{D}.
                               \end{cases}
\]
Adding the above two equations and using the boundary condition
\eqref{gflhe-bc}, we have
 \begin{equation}\label{gfl-if}
 \phi_{\rm inc}({\bf r})=\int_\Gamma G_\kappa({\bf r}, {\bf
r}')\partial_{\nu'}\phi({\bf r}'){\rm d}s({\bf r}')-\sum_{k=1}^m \sigma_k\phi_k
G_\kappa({\bf r}, {\bf r}_k),\quad {\bf r}\in D.
\end{equation}
and
\begin{equation}\label{gfl-sf}
\psi({\bf r})=\sum_{k=1}^m \sigma_k \phi_k G({\bf r}, {\bf r}_k)-\int_\Gamma
G_\kappa({\bf r}, {\bf r}')\partial_{\nu'}\phi({\bf r}'){\rm d}s({\bf r}'),\quad
{\bf r}\in\mathbb{R}^2\setminus\bar{D}.
\end{equation}
To compute the scattered field $\psi$, it is required to compute
$\partial_\nu\phi$ and $\phi_k, k=1, \dots, m$.

Adding the incident field on both sides of \eqref{gfl-sf} yields
\begin{equation}\label{gfl-tf}
\phi({\bf r})=\phi_{\rm inc}({\bf r})+\sum_{k=1}^m \sigma_k \phi_k G({\bf r},
{\bf r}_k)-\int_\Gamma G({\bf r}, {\bf r}')\partial_{\nu'}\phi({\bf
r}'){\rm d}s({\bf r}'),\quad {\bf r}\in\mathbb{R}^2\setminus\bar{D}.
\end{equation}
Evaluating \eqref{gfl-tf} at ${\bf r}_i$ and excluding the self-interaction of
the point scatterers, we get 
\begin{equation}\label{gfl-lps}
\phi_i-\sum_{\stackrel{\scriptstyle k=1}{k\ne
i}}^{m} \sigma_k \phi_k G_\kappa({\bf r}_i, {\bf r}_k)+\int_\Gamma G_\kappa({\bf
r}_i, {\bf r}')\partial_{\nu'}\phi({\bf r}'){\rm d}s({\bf r}')=\phi_{\rm
inc}({\bf r}_i).
\end{equation}
Taking the normal derivative of \eqref{gfl-if} on $\Gamma$ and applying the
jump condition lead to 
\begin{equation}\label{gfl-nif}
 \partial_\nu\phi_{\rm inc}({\bf r})=\int_\Gamma \partial_\nu G_\kappa({\bf r},
{\bf r}')\partial_{\nu'}\phi({\bf r}'){\rm d}s({\bf
r}')-\sum_{k=1}^m \sigma_k\phi_k \partial_\nu G_\kappa({\bf r}, {\bf
r}_k)+\frac{1}{2}\partial_\nu\phi({\bf r}),\quad {\bf r}\in\Gamma.
\end{equation}
Multiplying \eqref{gfl-if} by ${\rm i}\eta$ and subtracting it from
\eqref{gfl-nif}, we obtain 
\begin{equation}\label{gfl-li}
 \frac{1}{2}\partial_\nu\phi({\bf r})+\int_\Gamma (\partial_\nu-{\rm
i}\eta) G_\kappa({\bf r}, {\bf r}')\partial_{\nu'}\phi({\bf r}'){\rm
d}s({\bf r}')-\sum_{k=1}^m \sigma_k\phi_k (\partial_\nu-{\rm i}\eta)
G_\kappa({\bf r}, {\bf r}_k)=(\partial_\nu-{\rm i}\eta)\phi_{\rm inc}({\bf r}).
\end{equation}
The coupled system \eqref{gfl-lps} and \eqref{gfl-li} forms the generalized
Foldy--Lax formulation for the scattering problem with the linear point
scatterers and extended scatterers. 

\subsection{Quadratically nonlinear point scatterers}

Consider the point scatterers with weak quadratic nonlinearity. Let
$\kappa_j=\omega_j/c$ and $\phi^{(j)}$ be the field corresponding to the
wavenumber $\kappa_j$. Viewing the external field acting on the nonlinear point
scatterers as point sources for the extended obstacle, we have the following
equations in the exterior domain $\mathbb R^2\setminus\bar{D}$:
\begin{subequations}\label{gflhe-qtf}
 \begin{align}
  \Delta \phi^{(1)}({\bf r})+\kappa_1^2 \phi^{(1)}({\bf r})&=-\sum_{k=1}^m
\left(\sigma_{k,1}^{(1)}\phi_k^{(1)}+\sigma_{k,1}^{(2)}\bar{\phi}_k^{(1)}\phi_k^{(2)}
\right)\delta( { \bf r}-{\bf r}_k),\\
 \Delta \phi^{(2)}({\bf r})+\kappa_2^2 \phi^{(2)}({\bf r})&=-\sum_{k=1}^m
\left(\sigma_{k,2}^{(1)}\phi_k^{(2)}+\sigma_{k,2}^{(2)}\bigl(\phi_k^{(1)}\bigr)^2
\right)\delta({ \bf r}-{\bf r}_k).
 \end{align}
\end{subequations}
The sound-soft boundary condition implies that 
\begin{equation}
 \phi^{(1)}=\phi^{(2)}=0\quad\text{on}~\Gamma. 
\end{equation}
The fields satisfy the following relationship: 
\[
 \phi^{(1)}=\phi_{\rm inc}+\psi^{(1)},\quad \phi^{(2)}=\psi^{(2)},
\]
where $\psi^{(j)}$ is the scattered field corresponding to the wavenumber
$\kappa_j$ and satisfies the Sommerfeld radiation condition
\[
 \lim_{r\to\infty}r^{1/2}\left(\partial_r\psi^{(j)}-{\rm
i}\kappa_j\psi^{(j)}\right)=0,\quad r=|{\bf r}|.
\]

Similarly, we may show that the incident field satisfies
\[
 \int_\Gamma \left(\phi_{\rm inc}({\bf r}')\partial_{\nu'} G_{\kappa_1}({\bf r},
{\bf r}')-\partial_{\nu'}\phi_{\rm inc}({\bf r}')G_{\kappa_1}({\bf r}, {\bf
r}')\right) {\rm d}s({\bf r}')=\begin{cases}
                                -\phi_{\rm inc}({\bf r}),&\quad{\bf r}\in D,\\
                                0,&\quad{\bf r}\in\mathbb R^2\setminus\bar{D};
                               \end{cases}
\]
the scattered fields satisfy 
\begin{align*}
&\int_\Gamma \left(\psi^{(1)}({\bf r}')\partial_{\nu'} G_{\kappa_1}({\bf r},
{\bf r}')-\partial_{\nu'}\psi^{(1)}({\bf r}')G_{\kappa_1}({\bf r}, {\bf
r}')\right) {\rm d}s({\bf r}')\\
&\qquad+
\sum_{k=1}^m\left(
\sigma_{k,1}^{(1)}\phi_k^{(1)}+\sigma_{k,1}^{(2)}\bar{\phi}_k^ {(1)}\phi_k^{(2)} 
\right) G_{\kappa_1}({\bf r}, {\bf r}_k)=\begin{cases}
                                0,&\quad{\bf r}\in D,\\
                                \psi^{(1)}({\bf r}),&\quad{\bf r}\in\mathbb
R^2\setminus\bar{D}.
                               \end{cases}
\end{align*}
and
\begin{align*}
&\int_\Gamma \left(\psi^{(2)}({\bf r}')\partial_{\nu'} G_{\kappa_2}({\bf r},
{\bf r}')-\partial_{\nu'}\psi^{(2)}({\bf r}')G_{\kappa_2}({\bf r}, {\bf
r}')\right) {\rm d}s({\bf r}')\\
&\qquad+ \sum_{k=1}^m
\left(\sigma_{k,2}^{(1)}\phi_k^{(2)}+\sigma_{k,2}^{(2)}\bigl(\phi_k^{(1)}\bigr)^2
\right)G_{\kappa_2}({\bf r}, {\bf r}_k)=\begin{cases}
                                0,&\quad{\bf r}\in D,\\
                                \psi^{(2)}({\bf r}),&\quad{\bf r}\in\mathbb
R^2\setminus\bar{D}.
                               \end{cases}
\end{align*}
Adding the above equations and using the boundary condition yields
\begin{subequations}\label{gfl-qif}
\begin{align}
 \phi_{\rm inc}({\bf r})&=\int_\Gamma G_{\kappa_1}({\bf r}, {\bf
r}')\partial_{\nu'}\phi^{(1)}({\bf r}'){\rm d}s({\bf
r}')-\sum_{k=1}^m\left(
\sigma_{k,1}^{(1)}\phi_k^{(1)}+\sigma_{k,1}^{(2)}\bar{\phi}_k^{(1)}\phi_k^{(2)} 
\right) G_{\kappa_1}({\bf r}, {\bf r}_k),\quad {\bf r}\in
D,\\
0&=\int_\Gamma G_{\kappa_2}({\bf r}, {\bf
r}'){\rm d}s({\bf r}')\partial_{\nu'}\phi^{(2)}({\bf r}')- \sum_{k=1}^m
\left(\sigma_{k,2}^{(1)}\phi_k^{(2)}+\sigma_{k,2}^{(2)}\bigl(\phi_k^{(1)}\bigr)^2
\right)G_{\kappa_2}({\bf r}, {\bf r}_k) ,\quad {\bf r}\in
D,
\end{align}
\end{subequations}
and
\begin{subequations}\label{gfl-qsf}
 \begin{align}
\label{gfl-qsf1}\psi^{(1)}({\bf
r})&=\sum_{k=1}^m\left(
\sigma_{k,1}^{(1)}\phi_k^{(1)}+\sigma_{k,1}^{(2)}\bar{\phi}_k^{(1)}\phi_k^{(2)} 
\right) G_{\kappa_1}({\bf r}, {\bf r}_k)\notag\\
&\qquad-\int_\Gamma G_{\kappa_1}({\bf r}, {\bf
r}')\partial_{\nu'}\phi^{(1)}({\bf r}'){\rm d}s({\bf r}'),\quad {\bf
r}\in\mathbb{R}^2\setminus\bar{D},\\
\label{gfl-qsf2}\psi^{(2)}({\bf r})&=\sum_{k=1}^m
\left(\sigma_{k,2}^{(1)}\phi_k^{(2)}+\sigma_{k,2}^{(2)}\bigl(\phi_k^{(1)}\bigr)^2
\right)G_{\kappa_2}({\bf r}, {\bf r}_k)\notag\\
&\qquad-\int_\Gamma G_{\kappa_2}({\bf r}, {\bf
r}')\partial_{\nu'}\phi^{(2)}({\bf
r}'){\rm d}s({\bf r}'),\quad {\bf r}\in\mathbb{R}^2\setminus\bar{D}.
 \end{align}
\end{subequations}

Adding the incident field to \eqref{gfl-qsf1} and noting
$\phi^{(2)}=\psi^{(2)}$, we obtain 
\begin{subequations}\label{gfl-qtf}
 \begin{align}
\phi^{(1)}({\bf
r})&=\phi_{\rm
inc}({\bf r})+\sum_{k=1}^m\left(
\sigma_{k,1}^{(1)}\phi_k^{(1)}+\sigma_{k,1}^{(2)}\bar{\phi}_k^{(1)}\phi_k^{(2)} 
\right) G_{\kappa_1}({\bf r}, {\bf r}_k)\notag\\
&\qquad-\int_\Gamma G_{\kappa_1}({\bf r}, {\bf
r}')\partial_{\nu'}\phi^{(1)}({\bf r}'){\rm d}s({\bf r}'),\quad {\bf
r}\in\mathbb{R}^2\setminus\bar{D},\\
\phi^{(2)}({\bf r})&=\sum_{k=1}^m
\left(\sigma_{k,2}^{(1)}\phi_k^{(2)}+\sigma_{k,2}^{(2)}\bigl(\phi_k^{(1)}\bigr)^2
\right)G_{\kappa_2}({\bf r}, {\bf r}_k)\notag\\
&\qquad-\int_\Gamma G_{\kappa_2}({\bf r}, {\bf
r}')\partial_{\nu'}\phi^{(2)}({\bf
r}'){\rm d}s({\bf r}'),\quad {\bf r}\in\mathbb{R}^2\setminus\bar{D}.
 \end{align}
\end{subequations}
Evaluating \eqref{gfl-qtf} at ${\bf r}_i$ leads to 
\begin{subequations}\label{gfl-qps}
 \begin{align}
\phi^{(1)}_i&-\sum_{\stackrel{\scriptstyle k=1}{k\ne
i}}^{m}\left(
\sigma_{k,1}^{(1)}\phi_k^{(1)}+\sigma_{k,1}^{(2)}\bar{\phi}_k^{(1)}\phi_k^{(2)} 
\right) G_{\kappa_1}({\bf r}_i, {\bf
r}_k)\notag\\
&\qquad+\int_\Gamma G_{\kappa_1}({\bf r}_i, {\bf
r}')\partial_{\nu'}\phi^{(1)}({\bf r}'){\rm d}s({\bf r}')=\phi_{\rm
inc}({\bf r}_i),\\
\phi^{(2)}_i&-\sum_{\stackrel{\scriptstyle k=1}{k\ne
i}}^{m} \left(
\sigma_{k,2}^{(1)}\phi_k^{(2)}+\sigma_{k,2}^{(2)}\bigl(\phi_k^{(1)}\bigr)^2
\right)G_{\kappa_2}({\bf r}_i, {\bf r}_k)\notag\\
&\qquad+\int_\Gamma G_{\kappa_2}({\bf r}_i, {\bf
r}')\partial_{\nu'}\phi^{(2)}({\bf r}'){\rm d}s({\bf r}')=0. 
 \end{align}
\end{subequations}
Taking the normal derivative of \eqref{gfl-qif} and using the jump conditions,
we get 
\begin{subequations}\label{gfl-qnif}
\begin{align}
 \partial_\nu\phi_{\rm inc}({\bf r})=&\int_\Gamma \partial_\nu G_{\kappa_1}({\bf
r}, {\bf r}')\partial_{\nu'}\phi^{(1)}({\bf r}'){\rm d}s({\bf
r}')\notag\\
&\qquad-\sum_{k=1}^m\left(
\sigma_{k,1}^{(1)}\phi_k^{(1)}+\sigma_{k,1}^{(2)}\bar{\phi}_k^{(1)}\phi_k^{(2)} 
\right)\partial_\nu G_{\kappa_1}({\bf r}, {\bf
r}_k)+\frac{1}{2}\partial_\nu\phi^{(1)}({\bf r}),\\
0=&\int_\Gamma \partial_\nu G_{\kappa_2}({\bf
r}, {\bf r}')\partial_{\nu'}\phi^{(2)}({\bf r}'){\rm d}s({\bf
r}')\notag\\
&\qquad- \sum_{k=1}^m
\left(\sigma_{k,2}^{(1)}\phi_k^{(2)}+\sigma_{k,2}^{(2)}\bigl(\phi_k^{(1)}\bigr)^2
\right)\partial_\nu G_{\kappa_2}({\bf r}, {\bf
r}_k)+\frac{1}{2}\partial_\nu\phi^{(2)}({\bf r}).
\end{align}
\end{subequations}
Multiplying \eqref{gfl-qif} by ${\rm i}\eta$ and subtract it from
\eqref{gfl-qnif} give
\begin{subequations}\label{gfl-qi}
\begin{align}
\frac{1}{2}\partial_\nu\phi^{(1)}({\bf r})&+\int_\Gamma (\partial_\nu-{\rm
i}\eta) G_{\kappa_1}({\bf r}, {\bf r}')\partial_{\nu'}\phi^{(1)}({\bf r}'){\rm
d}s({\bf r}')\notag\\
&\qquad-\sum_{k=1}^m\left(
\sigma_{k,1}^{(1)}\phi_k^{(1)}+\sigma_{k,1}^{(2)}\bar{\phi}_k^ {(1)} \phi_k^{(2)} 
\right)(\partial_\nu-{\rm i}\eta) G_{\kappa_1}({\bf
r}, {\bf r}_k)=(\partial_\nu-{\rm i}\eta)\phi_{\rm inc}({\bf r}),\\
\frac{1}{2}\partial_\nu\phi^{(2)}({\bf r})&+\int_\Gamma (\partial_\nu-{\rm
i}\eta) G_{\kappa_2}({\bf r}, {\bf r}')\partial_{\nu'}\phi^{(2)}({\bf r}'){\rm
d}s({\bf r}') \notag\\
&\qquad - \sum_{k=1}^m
\left(\sigma_{k,2}^{(1)}\phi_k^{(2)}+\sigma_{k,2}^{(2)}\bigl(\phi_k^{(1)}\bigr)^2
\right)(\partial_\nu -{\rm i}\eta)G_{\kappa_2}({\bf r}, {\bf r}_k)=0.
\end{align}
\end{subequations}
The coupled system \eqref{gfl-qps} and \eqref{gfl-qi} gives the generalized
Foldy--Lax formulation for the scattering problem with quadratic nonlinear
point scatterers and extended scatterers. 

\subsection{Cubically nonlinear point scatterers}

Consider the point scatterers with weak cubic nonlinearity. Let
$\kappa_j=\omega_j/c$ and $\phi^{(j)}$ be the field corresponding to the
wavenumber $\kappa_j$. Viewing the external field acting on the nonlinear point
scatterers as point sources for the extended obstacle, we have the following
equations in the exterior domain $\mathbb R^2\setminus\bar{D}$:
\begin{subequations}\label{gflhe-ctf}
 \begin{align}
  \Delta \phi^{(1)}({\bf r})+\kappa_1^2 \phi^{(1)}({\bf r})&=-\sum_{k=1}^m
\left(\sigma_{k,1}^{(1)}\phi^{(1)}_k+\sigma_{k,1}^{(3)}|\phi_k^{(1)}|^2\phi_k^{(1)}+\sigma_{k,2}^{(3)}\bigl(\bar{\phi}_k^{(1)}\bigr)^2\phi_k^{(3)}
\right)\delta({ \bf r}-{\bf r}_k),\\
 \Delta \phi^{(3)}({\bf r})+\kappa_3^2 \phi^{(3)}({\bf r})&=-\sum_{k=1}^m
\left(\sigma_{k,2}^{(1)}\phi_k^{(3)}+\sigma_{k,3}^{(3)}\bigl(\phi_k^{(1)}\bigr)^3
\right)\delta({\bf r}-{\bf r}_k).
 \end{align}
\end{subequations}
The sound-soft boundary condition implies that 
\begin{equation}
 \phi^{(1)}=\phi^{(3)}=0\quad\text{on}~\Gamma. 
\end{equation}
The fields satisfy the following relationship: 
\[
 \phi^{(1)}=\phi_{\rm inc}+\psi^{(1)},\quad \phi^{(3)}=\psi^{(3)},
\]
where $\psi^{(j)}$ is the scattered field corresponding to the wavenumber
$\kappa_j$ and satisfies the Sommerfeld radiation condition
\[
 \lim_{r\to\infty}r^{1/2}\left(\partial_r\psi^{(j)}-{\rm
i}\kappa_j\psi^{(j)}\right)=0,\quad r=|{\bf r}|.
\]

Similarly, we may show that the incident field satisfies
\[
 \int_\Gamma \left(\phi_{\rm inc}({\bf r}')\partial_{\nu'} G_{\kappa_1}({\bf r},
{\bf r}')-\partial_{\nu'}\phi_{\rm inc}({\bf r}')G_{\kappa_1}({\bf r}, {\bf
r}')\right) {\rm d}s({\bf r}')=\begin{cases}
                                -\phi_{\rm inc}({\bf r}),&\quad{\bf r}\in D,\\
                                0,&\quad{\bf r}\in\mathbb R^2\setminus\bar{D};
                               \end{cases}
\]
the scattered fields satisfy 
\begin{align*}
&\int_\Gamma \left(\psi^{(1)}({\bf r}')\partial_{\nu'} G_{\kappa_1}({\bf r},
{\bf r}')-\partial_{\nu'}\psi^{(1)}({\bf r}')G_{\kappa_1}({\bf r}, {\bf
r}')\right) {\rm d}s({\bf r}')\\
&\qquad+
\sum_{k=1}^m
\left(\sigma_{k,1}^{(1)}\phi^{(1)}_k+\sigma_{k,1}^{(3)}|\phi_k^{(1)}|^2\phi_k^{(1)}+\sigma_{k,2}^{(3)}\bigl(\bar{\phi}_k^{(1)}\bigr)^2\phi_k^{(3)}
\right) G_{\kappa_1}({\bf r}, {\bf r}_k)=\begin{cases}
                                0,&\quad{\bf r}\in D,\\
                                \psi^{(1)}({\bf r}),&\quad{\bf r}\in\mathbb
R^2\setminus\bar{D}.
                               \end{cases}
\end{align*}
and
\begin{align*}
&\int_\Gamma \left(\psi^{(3)}({\bf r}')\partial_{\nu'} G_{\kappa_3}({\bf r},
{\bf r}')-\partial_{\nu'}\psi^{(3)}({\bf r}')G_{\kappa_3}({\bf r}, {\bf
r}')\right) {\rm d}s({\bf r}')\\
&\qquad+ \sum_{k=1}^m
\left(\sigma_{k,2}^{(1)}\phi_k^{(3)}+\sigma_{k,3}^{(3)}\bigl(\phi_k^{(1)}\bigr)^3
\right) G_{\kappa_3}({\bf r}, {\bf r}_k)=\begin{cases}
                                0,&\quad{\bf r}\in D,\\
                                \psi^{(3)}({\bf r}),&\quad{\bf r}\in\mathbb
R^2\setminus\bar{D}.
                               \end{cases}
\end{align*}
Adding the above equations and using the boundary condition yield
\begin{subequations}\label{gfl-cif}
\begin{align}
 \phi_{\rm inc}({\bf r})&=\int_\Gamma G_{\kappa_1}({\bf r}, {\bf
r}')\partial_{\nu'}\phi^{(1)}({\bf r}'){\rm d}s({\bf
r}')\notag\\
&\qquad-\sum_{k=1}^m
\left(\sigma_{k,1}^{(1)}\phi^{(1)}_k+\sigma_{k,1}^{(3)}|\phi_k^{(1)}|^2\phi_k^{(1)}+\sigma_{k,2}^{(3)}\bigl(\bar{\phi}_k^{(1)}\bigr)^2\phi_k^{(3)}
\right) G_{\kappa_1}({\bf r}, {\bf r}_k),\quad {\bf r}\in
D,\\
0&=\int_\Gamma G_{\kappa_3}({\bf r}, {\bf
r}'){\rm d}s({\bf r}')\partial_{\nu'}\phi^{(3)}({\bf r}')- \sum_{k=1}^m
\left(\sigma_{k,2}^{(1)}\phi_k^{(3)}+\sigma_{k,3}^{(3)}\bigl(\phi_k^{(1)}\bigr)^3
\right)G_{\kappa_3}({\bf r}, {\bf r}_k) ,\quad {\bf r}\in
D,
\end{align}
\end{subequations}
and
\begin{subequations}\label{gfl-csf}
 \begin{align}
\label{gfl-csf1}\psi^{(1)}({\bf
r})&=\sum_{k=1}^m
\left(\sigma_{k,1}^{(1)}\phi^{(1)}_k+\sigma_{k,1}^{(3)}|\phi_k^{(1)}|^2\phi_k^{(1)}+\sigma_{k,2}^{(3)}\bigl(\bar{\phi}_k^{(1)}\bigr)^2\phi_k^{(3)}
\right) G_{\kappa_1}({\bf r}, {\bf r}_k)\notag\\
&\qquad-\int_\Gamma G_{\kappa_1}({\bf r}, {\bf
r}')\partial_{\nu'}\phi^{(1)}({\bf r}'){\rm d}s({\bf r}'),\quad {\bf
r}\in\mathbb{R}^2\setminus\bar{D},\\
\label{gfl-csf2}\psi^{(3)}({\bf r})&=\sum_{k=1}^m
\left(\sigma_{k,2}^{(1)}\phi_k^{(3)}+\sigma_{k,3}^{(3)}\bigl(\phi_k^{(1)}\bigr)^3
\right)G_{\kappa_3}({\bf r}, {\bf r}_k)\notag\\
&\qquad-\int_\Gamma G_{\kappa_3}({\bf r}, {\bf
r}')\partial_{\nu'}\phi^{(3)}({\bf
r}'){\rm d}s({\bf r}'),\quad {\bf r}\in\mathbb{R}^2\setminus\bar{D}.
 \end{align}
\end{subequations}

Adding the incident field to \eqref{gfl-csf1} and noting
$\phi^{(3)}=\psi^{(3)}$, we obtain 
\begin{subequations}\label{gfl-ctf}
 \begin{align}
\phi^{(1)}({\bf
r})&=\phi_{\rm
inc}({\bf r})+\sum_{k=1}^m
\left(\sigma_{k,1}^{(1)}\phi^{(1)}_k+\sigma_{k,1}^{(3)}|\phi_k^{(1)}|^2\phi_k^{(1)}+\sigma_{k,2}^{(3)}\bigl(\bar{\phi}_k^{(1)}\bigr)^2\phi_k^{(3)}
\right) G_{\kappa_1}({\bf r}, {\bf r}_k)\notag\\
&\qquad-\int_\Gamma G_{\kappa_1}({\bf r}, {\bf
r}')\partial_{\nu'}\phi^{(1)}({\bf r}'){\rm d}s({\bf r}'),\quad {\bf
r}\in\mathbb{R}^2\setminus\bar{D},\\
\phi^{(3)}({\bf r})&=\sum_{k=1}^m
\left(\sigma_{k,2}^{(1)}\phi_k^{(3)}+\sigma_{k,3}^{(3)}\bigl(\phi_k^{(1)}\bigr)^3
\right) G_{\kappa_3}({\bf r}, {\bf
r}_k)\notag\\
&\qquad-\int_\Gamma G_{\kappa_3}({\bf r}, {\bf
r}')\partial_{\nu'}\phi^{(3)}({\bf
r}'){\rm d}s({\bf r}'),\quad {\bf r}\in\mathbb{R}^2\setminus\bar{D}.
 \end{align}
\end{subequations}
Evaluating \eqref{gfl-ctf} at ${\bf r}_i$ and excluding the
self-interaction lead to 
\begin{subequations}\label{gfl-cps}
 \begin{align}
\phi^{(1)}_i&-\sum_{\stackrel{\scriptstyle k=1}{k\ne
i}}^{m}\left(\sigma_{k,1}^{(1)}\phi_k^{(1)}+\sigma_{k,1}^{(3)}|\phi_k^{(1)}|^2\phi_k^{(1)} +\sigma_{k,2}^{(3)}\bigl(\bar{\phi}_k^{(1)} \bigr)^2\phi_k^{(3)}
\right) G_{\kappa_1}({\bf r}_i, {\bf
r}_k)\notag\\
&\qquad+\int_\Gamma G_{\kappa_1}({\bf r}_i, {\bf
r}')\partial_{\nu'}\phi^{(1)}({\bf r}'){\rm d}s({\bf r}')=\phi_{\rm
inc}({\bf r}_i),\\
\phi^{(3)}_i&-\sum_{\stackrel{\scriptstyle k=1}{k\ne
i}}^{m}\left(\sigma_{k,2}^{(1)}\phi_k^{(3)}+\sigma_{k,3}^{(3)}\bigl(\phi_k^{(1)}\bigr)^3
\right) G_{\kappa_3}({\bf r}_i, {\bf
r}_k)\notag\\
&\qquad+\int_\Gamma G_{\kappa_3}({\bf r}_i, {\bf
r}')\partial_{\nu'}\phi^{(3)}({\bf r}'){\rm d}s({\bf r}')=0. 
 \end{align}
\end{subequations}
Taking the normal derivative of \eqref{gfl-cif} and using the jump conditions,
we get 
\begin{subequations}\label{gfl-cnif}
\begin{align}
 \partial_\nu\phi_{\rm inc}({\bf r})&=\int_\Gamma \partial_\nu G_{\kappa_1}({\bf
r}, {\bf r}')\partial_{\nu'}\phi^{(1)}({\bf r}'){\rm d}s({\bf
r}')\notag\\
&\quad-\sum_{k=1}^m
\left(\sigma_{k,1}^{(1)}\phi^{(1)}_k+\sigma_{k,1}^{(3)}|\phi_k^{(1)}|^2\phi_k^{(1)}+\sigma_{k,2}^{(3)}\bigl(\bar{\phi}_k^{(1)}\bigr)^2\phi_k^{(3)}
\right) \partial_\nu G_{\kappa_1}({\bf r}, {\bf
r}_k)+\frac{1}{2}\partial_\nu\phi^{(1)}({\bf r}),\\
0=&\int_\Gamma \partial_\nu G_{\kappa_3}({\bf
r}, {\bf r}')\partial_{\nu'}\phi^{(3)}({\bf r}'){\rm d}s({\bf
r}')\notag\\
&\qquad- \sum_{k=1}^m
\left(\sigma_{k,2}^{(1)}\phi_k^{(3)}+\sigma_{k,3}^{(3)}\bigl(\phi_k^{(1)}\bigr)^3
\right)\partial_\nu G_{\kappa_3}({\bf r}, {\bf
r}_k)+\frac{1}{2}\partial_\nu\phi^{(3)}({\bf r}).
\end{align}
\end{subequations}
Multiplying \eqref{gfl-qif} by ${\rm i}\eta$ and subtract it from
\eqref{gfl-qnif} give
\begin{subequations}\label{gfl-ci}
\begin{align}
\frac{1}{2}\partial_\nu\phi^{(1)}({\bf
r})&-\sum_{k=1}^m\left(\sigma_{k,1}^{(1)}\phi^{(1)}_k+\sigma_{k,1}^{(3)}|\phi_k^{(1)}|^2\phi_k^{(1)}+\sigma_{k,2}^{(3)}\bigl(\bar{\phi}_k^{(1)}\bigr)^2\phi_k^{(3)}
\right)(\partial_\nu-{\rm i}\eta) G_{\kappa_1}({\bf
r}, {\bf r}_k)\notag\\
&\quad+\int_\Gamma (\partial_\nu-{\rm
i}\eta) G_{\kappa_1}({\bf r}, {\bf r}')\partial_{\nu'}\phi^{(1)}({\bf r}'){\rm
d}s({\bf r}')=(\partial_\nu-{\rm i}\eta)\phi_{\rm inc}({\bf r}),\\
\frac{1}{2}\partial_\nu\phi^{(3)}({\bf r})&- \sum_{k=1}^m
\left(\sigma_{k,2}^{(1)}\phi_k^{(3)}+\sigma_{k,3}^{(3)}\bigl(\phi_k^{(1)}\bigr)^3
\right)(\partial_\nu -{\rm i}\eta)G_{\kappa_3}({\bf r}, {\bf r}_k)\notag\\
&\quad +\int_\Gamma (\partial_\nu-{\rm
i}\eta) G_{\kappa_3}({\bf r}, {\bf r}')\partial_{\nu'}\phi^{(3)}({\bf r}'){\rm
d}s({\bf r}')=0.
\end{align}
\end{subequations}
The coupled system \eqref{gfl-cps} and \eqref{gfl-ci} gives the generalized
Foldy--Lax formulation for the scattering problem with cubic nonlinear
point scatterers and extended scatterers. 

\section{Direct imaging method}

In this section, we introduce a fast direct imaging method to reconstruct the
shape of the extended scatterers. 

\subsection{Far-field pattern}

The far-field pattern is an important quantity which encodes the information
about the scatterers such as location and shape. Given an incident field with
incident direction ${\bf d}$, the scattered field has the asymptotic behavior
\begin{equation}\label{ffp}
 \psi({\bf r}, {\bf d})=\frac{e^{{\rm i}\kappa|{\bf r}|}}{|{\bf
r}|^{\frac{1}{2}}}\left(\psi_{\infty}(\hat{\bf r}, {\bf d})+O(|{\bf
r}|)^{-1}\right)\quad\text{as}\,|{\bf r}|\to\infty,
\end{equation}
uniformly in all directions $\hat{\bf r}={\bf r}/|{\bf r}|$, where
the function $\psi_{\infty}$ is called the far-field pattern of the scattered
field $\psi$, and $\hat{\bf r}\in\mathbb{S}^1$ is the observation direction.

Recall the asymptotic behavior for the Hankel function for large arguments
\[
H_0^{(1)}(z)=\sqrt{\frac{2}{\pi z}}\,e^{{\rm
i} \left(z-\frac{\pi}{4}\right)}\left(1+O(z^{-1})\right)
\]
and the following identity
\[
|{\bf r}-{\bf r}'|=\sqrt{|{\bf r}|^2-2{\bf r}\cdot{\bf r}'+|{\bf
r}'|^2}=|{\bf r}|-\hat{\bf r}\cdot{\bf r}'+O(|{\bf
r}|^{-1})\quad\text{as}\,{|{\bf r}|}\to\infty.
\]
Using \eqref{ffp} and the scattered field representations \eqref{fl-lsf},
\eqref{fl-qsf}, \eqref{fl-csf}, \eqref{bi-sf}, \eqref{gfl-sf}, \eqref{gfl-qsf},
\eqref{gfl-csf}, we obtain the following far-field patterns of the
scattered field for the scattering problem with point scatterers, extended
scatterers, and mixed scatterers, respectively. 

\begin{enumerate}

\item[(i)] Foldy--Lax formulation for point scatterers:

\begin{enumerate}
 
 \item[(a)] Linear point scatterers
 
 \begin{equation}\label{ffp-fl-lps}
 \psi_{\infty, {\rm FL, l}}(\hat{\bf r}, {\bf
d})=\gamma\sum_{k=1}^m \sigma_k \phi_k({\bf d}) e^{-{\rm i}\kappa\hat{\bf
r}\cdot{\bf r}_k};
\end{equation}
 
 \item[(b)] Quadratically nonlinear point scatterers
 
 \begin{subequations}\label{ffp-fl-qps}
 \begin{align}
  \psi^{(1)}_{\infty, {\rm FL, q}}(\hat{\bf r}, {\bf
d})&=\gamma_1\sum_{k=1}^m\left(
\sigma_{k,1}^{(1)}\phi_k^{(1)}+\sigma_{k,1}^{(2)}\bar{\phi}_k^{(1)}\phi_k^{(2)} 
\right) e^{-{\rm i}\kappa\hat{\bf r}\cdot{\bf
r}_k},\\
 \psi^{(2)}_{\infty, {\rm FL, q}}(\hat{\bf r}, {\bf d})&=\gamma_2
\sum_{k=1}^m \left(
\sigma_{k,2}^{(1)}\phi_k^{(2)}+\sigma_{k,2}^{(2)}\bigl(\phi_k^{(1)}\bigr)^2 
\right) e^{-{\rm i}\kappa\hat{\bf r}\cdot{\bf r}_k};
 \end{align}
\end{subequations}

 \item[(c)] Cubically nonlinear point scatterers
 
 \begin{subequations}\label{ffp-fl-cps}
 \begin{align}
  \psi^{(1)}_{\infty, {\rm FL, c}}(\hat{\bf r}, {\bf
d})=&\gamma_1\sum_{k=1}^m
\left(\sigma_{k,1}^{(1)}\phi^{(1)}_k+\sigma_{k,1}^{(3)}|\phi_k^{(1)}|^2\phi_k^{(1)}+\sigma_{k,2}^{(3)}\bigl(\bar{\phi}_k^{(1)}\bigr)^2\phi_k^{(3)}
\right) e^{-{\rm i}\kappa\hat{\bf r}\cdot{\bf
r}_k},\\
 \psi^{(3)}_{\infty, {\rm FL, c}}(\hat{\bf r}, {\bf d})=&\gamma_3 \sum_{k=1}^m
\left(\sigma_{k,2}^{(1)}\phi_k^{(3)}+\sigma_{k,3}^{(3)}\bigl(\phi_k^{(1)}
\bigr)^3 \right) e^{-{\rm i}\kappa\hat{\bf r}\cdot{\bf
r}_k};
 \end{align}
\end{subequations}
 
\end{enumerate}

\item[(ii)] Boundary integral formulation for extended scatterers

\begin{equation}\label{ffp-bi}
 \psi_{\infty, {\rm BI}}(\hat{\bf r}, {\bf
d})=-\gamma\int_{\Gamma}\partial_{\nu'} \phi({\bf r'}; {\bf d})e^{-{\rm
i}\kappa\hat{\bf r}\cdot{\bf r}'}{\rm d}s({\bf
r}');
\end{equation}

\item[(iii)] Generalized Foldy--Lax formulation for mixed scatterers

\begin{enumerate}

\item[(a)] Linear point scatterers

\begin{equation}\label{ffp-gfl-lps}
 \psi_{\infty, {\rm GFL, l}}(\hat{\bf r}, {\bf
d})=\gamma\left[\sum_{k=1}^m \sigma_k \phi_k({\bf d}) e^{-{\rm i}\kappa\hat{\bf
r}\cdot{\bf r}_k}-\int_{\Gamma}\partial_{\nu'}
\phi({\bf r'}; {\bf d})e^{-{\rm i}\kappa\hat{\bf r}\cdot{\bf r}'}{\rm d}s({\bf
r}')\right];
\end{equation}

\item[(b)] Quadratically nonlinear point scatterers

\begin{subequations}\label{ffp-gfl-qps}
 \begin{align}
  \psi^{(1)}_{\infty, {\rm GFL, q}}(\hat{\bf r}, {\bf
d})&=\gamma_1\left[ \sum_{k=1}^m\left(
\sigma_{k,1}^{(1)}\phi_k^{(1)}+\sigma_{k,1}^{(2)}\bar{\phi}_k^{(1)}\phi_k^{(2)} 
\right) e^{-{\rm i}\kappa\hat{\bf r}\cdot{\bf
r}_k}-\int_{\Gamma}\partial_{\nu'} \phi^{(1)}({\bf r'}; {\bf d})e^{-{\rm
i}\kappa\hat{\bf r}\cdot{\bf r}'}{\rm d}s({\bf r}')\right],\\
 \psi^{(2)}_{\infty, {\rm GFL, q}}(\hat{\bf r}, {\bf d})&=\gamma_2\left[
\sum_{k=1}^m \left(
\sigma_{k,2}^{(1)}\phi_k^{(2)}+\sigma_{k,2}^{(2)}\bigl(\phi_k^{(1)}\bigr)^2 
\right) e^{-{\rm i}\kappa\hat{\bf r}\cdot{\bf
r}_k}-\int_{\Gamma}\partial_{\nu'} \phi^{(2)}({\bf r'}; {\bf d})e^{-{\rm
i}\kappa\hat{\bf r}\cdot{\bf r}'}{\rm d}s({\bf r}')\right],
 \end{align}
\end{subequations}

\item[(c)] Cubically nonlinear point scatterers

\begin{subequations}\label{ffp-gfl-cps}
 \begin{align}
  \psi^{(1)}_{\infty, {\rm GFL, c}}(\hat{\bf r}, {\bf
d})=&\gamma_1\Biggl[ \sum_{k=1}^m
\left(\sigma_{k,1}^{(1)}\phi^{(1)}_k+\sigma_{k,1}^{(3)}|\phi_k^{(1)}|^2\phi_k^{(1)}+\sigma_{k,2}^{(3)}\bigl(\bar{\phi}_k^{(1)}\bigr)^2\phi_k^{(3)}
\right) e^{-{\rm i}\kappa\hat{\bf r}\cdot{\bf
r}_k}\notag\\
&\hspace{6.5cm}-\int_{\Gamma}\partial_{\nu'} \phi^{(1)}({\bf r'}; {\bf
d})e^{-{\rm i}\kappa\hat{\bf r}\cdot{\bf r}'}{\rm d}s({\bf r}')\Biggr],\\
 \psi^{(3)}_{\infty, {\rm GFL, c}}(\hat{\bf r}, {\bf d})=&\gamma_3\left[
\sum_{k=1}^m \left(\sigma_{k,2}^{(1)}\phi_k^{(3)}+\sigma_{k,3}^{(3)}\bigl(\phi_k^{(1)}\bigr)^3 \right) e^{-{\rm i}\kappa\hat{\bf r}\cdot{\bf
r}_k}-\int_{\Gamma}\partial_{\nu'} \phi^{(3)}({\bf r'}; {\bf d})e^{-{\rm
i}\kappa\hat{\bf r}\cdot{\bf r}'}{\rm d}s({\bf r}')\right],
 \end{align}
\end{subequations}

\end{enumerate}

\end{enumerate}
where  
\[
\gamma=\frac{e^{{\rm i}\frac{\pi}{4}}}{\sqrt{8\pi\kappa}}, \quad
\gamma_j=\frac{e^{{\rm i}\frac{\pi}{4}}}{\sqrt{8\pi\kappa_j}}.
\]
When the observation directions and the number of point scatterers are large,
it is very slow to directly evaluate the far-field patterns
\eqref{ffp-fl-lps}--\eqref{ffp-gfl-cps}. In practice, the evaluation of the
far-field patterns are accelerated by the fast multipole method (FMM)\cite{GR1987}.  

\subsection{Imaging function}

Consider an array of transmitters that can send out plane incident waves and
record the far-field pattern of the scattered waves. Assume that we have a set
of incident plane waves with incident directions ${\bf d}_1, \dots, {\bf d}_N$
and the far-field patterns are recorded at observation direction $\hat{\bf r}_1,
\dots, \hat{\bf r}_M$, where ${\bf r}_i=(\cos\alpha_i, \sin\alpha_i)$ and
${\bf d}_j=(\cos\beta_j, \sin\beta_j), i=1, \dots, M, j=1, \dots, N$. Here
$\alpha_i$ is the observation angle and $\beta_j$ is the incident angle. These
measurement of the far-field patterns form an $M\times N$ response matrix
\begin{equation}\label{rm}
P_\kappa=\begin{bmatrix}
          \psi_\infty(\kappa; \hat{\bf r}_1, {\bf d}_1) & \cdots &
\psi_\infty(\kappa; \hat{\bf r}_1, {\bf d}_N)\\
          \vdots & \vdots & \vdots\\
           \psi_\infty(\kappa; \hat{\bf r}_M, {\bf d}_1) &\cdots & 
\psi_\infty(\kappa; \hat{\bf r}_M, {\bf d}_N)
         \end{bmatrix},
\end{equation}
where the far-field pattern $\psi_\infty$ represents any one of the far-field
patterns in \eqref{ffp-fl-lps}--\eqref{ffp-gfl-cps}.

Consider two unit vectors:
\[
{\bf u}_\kappa({\bf r})=\frac{1}{\sqrt{M}}(e^{{\rm i}\kappa{\bf
r}\cdot\hat{\bf r}_1}, \dots, e^{{\rm i}\kappa{\bf r}\cdot\hat{\bf r}_M})^\top
\]
and
\[
{\bf v}_\kappa({\bf r})=\frac{1}{\sqrt{N}}(e^{{\rm i}\kappa{\bf
r}\cdot{\bf d}_1}, \dots, e^{{\rm i}\kappa{\bf r}\cdot{\bf d}_N})^\top,
\]
which is the illumination vector with respect to the receivers and the
transmitters, respectively. Define an imaging function
\begin{equation}\label{imf}
 I_\kappa({\bf r})={\bf u}^\top_\kappa({\bf r}) P_\kappa{\bf v}_\kappa({\bf
r}).
\end{equation}
The direct imaging method is to evaluate the imaging function \eqref{imf} at any
given sampling point ${\bf r}\in\mathbb R^2$. 

Since the imaging function $I({\bf r})$ needs to be evaluated at every sampling
point and each evaluation requires the matrix-vector multiplication, the
computation is intensive. However, the imaging function can be written as 
\begin{align}\label{imf-fft}
I_\kappa({\bf r}) & = \frac{1}{\sqrt{MN}}\sum_{i=1}^{M}\sum_{j
=1}^{N} e^{{\rm i}\kappa{\bf r}\cdot (\hat{\bf r}_i+{\bf d}_j)}P^{(i,
j)}_\kappa\notag\\
&=\frac{1}{\sqrt{MN}}\sum_{i=1}^{M}\sum_{j = 1}^{N} e^{{\rm
i}\kappa\left[x(\cos{\alpha_i}+\cos{\beta_j})
+y(\sin{\alpha_i}+\sin{\beta_j})\right]}P_\kappa^{(i, j)}.
\end{align} 
It is clear to note from \eqref{imf-fft} that the imaging function
$I_\kappa({\bf r})$ is the two-dimensional Fourier transform of the $(i, j)$-th
entry of the response matrix $P_\kappa^{(i, j)}$, which takes the response
matrix from the frequency space $(\kappa(\cos{\alpha_i}+\cos{\beta_j}),
\kappa(\sin{\alpha_i}+\sin{\beta_j}))$ to the physical space $(x,y)$. While the
frequency points $(\kappa(\cos{\alpha_i}+\cos{\beta_j}),
\kappa(\sin{\alpha_i}+\sin{\beta_j}))$ are not uniform, one needs to apply the
two-dimensional non-uniform fast Fourier transform (NUFFT) to accelerate the
evaluation. 

\subsection{NUFFT}

Here we give a short introduction to NUFFT in one dimension. Details can be
found in \cite{gl-sr04}. Application of NUFFT on the acceleration of evaluating
forward multi-particle scattering in a layered medium can be found
in\cite{Lai:14}.  

Consider the following expression:
\begin{equation}\label{fi}
f(x_i) = \sum_{j=1}^{n}c_j e^{{\rm i}\xi_j x_i}, \quad i = 1,\dots, m,
\end{equation}
where $c_j\in\mathbb C, x_i\in [-m/2,m/2-1], \xi_j\in[-\pi,\pi]$. Given $c_j$ and
$\xi_j$, the goal is to evaluate $f(x_i)$ efficiently. When $m=n$ and $x_i$,
$\xi_j$ are uniformly distributed in their intervals, the sum can be accelerated
by standard FFT.

Consider the situation where $m\neq n$ and $\xi_j$ is not uniformly given in the
interval $[-\pi,\pi]$. Assume that $x_i$ is still uniformly distributed. We first
rewrite equation \eqref{fi} as
\[
f(x_i) =  \sum_{j=1}^{n}c_j \int_{-\infty}^{\infty} \delta(\xi-\xi_j)e^{{\rm i}
\xi x_i}{\rm d}\xi = \int_{-\infty}^{\infty} F(\xi) e^{{\rm i}\xi x_i}{\rm d}\xi,
\quad i = 1, \dots, m,
\]
where $\delta$ is the Dirac delta function and 
\[
F(\xi) = \sum_{j=1}^{n}c_j \delta(\xi-\xi_j).
\]
In other words, $f(x)$ can be taken as the Fourier transform of $F(\xi)$.
However, $F(\xi)$ consists of delta functions, which is
numerically difficult to evaluate. The function $F(\xi)$ can be modified by
convolving it with Gaussian function $g(\xi) = e^{-\xi^2/{4\tau}}$, where
$\tau>0$ is a small number. Let 
\[
\tilde{F}(\xi) =  \int_{-\infty}^{\infty}F(t)e^{-(\xi-t)^2/{4\tau}}{\rm d}t =
\sum_{j=1}^{n}c_j e^{-(\xi-\xi_j)^2/{4\tau}}. 
\]
We over sample $\tilde{F}(\xi)$ uniformly on the interval $[-\pi,\pi]$ and take
the Fourier transform of $\tilde{F}(\xi)$ with the standard FFT. Notice that to
apply FFT, both $F(\xi)$ and $g(\xi)$ are periodized first on $[-\pi,\pi]$ 
before evaluating $\tilde{F}(\xi)$\cite{gl-sr04}. Denote by 
$\hat{\tilde{F}}(x)$ and $\hat{g}(x)$ the Fourier transform
$\tilde{F}(\xi)$ and $g(\xi)$, respectively. It follows from the convolution
property of Fourier transform that $f(x)$ can be approximated by 
\[
f(x)=\frac{\hat{\tilde{F}}(x)}{\hat{g}(x)}=\left(\frac{\tau}{\pi}\right)^{
-\frac{1}{2}}e^{\tau x^2} \hat{\tilde F}(x).
\]

To conclude, the NUFFT consists of three steps: (1) convolution with Gaussian
function, (2) standard FFT for the over sampled function, (3) deconvolution with Gaussian
function. The choices of parameter $\tau$ and oversampling factor in step 2
highly affect the numerical performance of the NUFFT. For the one dimensional
case, $\tau$ is chosen to be $12/m^2$ and the number of over sampling points is $2m$, which
guarantees 12 digits precision. Using these parameters, we can easily find out
that the complexity of the one-dimensional NUFFT is $O(K\log K)$, where
$K=\max\{m,n\}$. The scheme introduced here can be easily extended to the two-
and higher-dimensional cases. The two-dimensional NUFFT is used in this work to
evaluate \eqref{imf-fft}.

\subsection{Imaging with nonlinear point scatterers}

\begin{figure}
\centering
\includegraphics[width=0.35\textwidth]{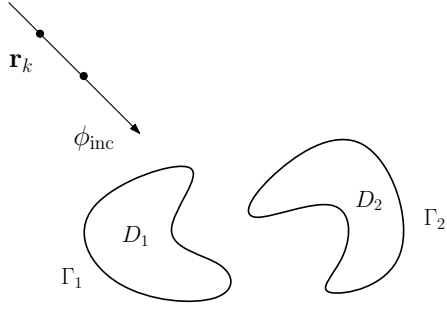}
\caption{Schematic of the imaging modality with nonlinear point scatterers.}
\label{pg_np}
\end{figure}

The method introduced in the beginning of this section works effectively to
image the extended scatterers surrounded by linear point scatterers. The method
can also be applied directly to the scattered field generated by the first order
frequency wave in the nonlinear case. The purpose of adding nonlinear point
scatterers is to introduce higher order field, which allows to capture more
details of the extended scatterers. 

However, if we fix the location of the nonlinear point scatterers, the imaging
method will fail for testing the higher order frequency waves. The reason is
obvious: in the direct imaging method, the test function ${\bf v}_\kappa$
depends on the incident angle ${\bf d}_j, j=1, \dots, N$. For higher order
frequency waves, the incident wave is essentially generated by the nonlinear
interaction of point scatterers. The direction of such generated incident wave
does not line up with the incident direction in the test function.

To remedy the scheme, we first put the point scatterers sufficiently far away
from the extended scatterers and line up the point scatterers towards the
incident direction, as is shown in Figure \ref{pg_np}. In addition, we move the
point scatterers along with the change of incident direction. For the
higher order frequency waves, the entries of the response matrix are taken as
the difference of the far-field patterns from the generalized Foldy--Lax
formulation and the Foldy--Lax formulation, i.e., 
\[
 P_\kappa^{(i, j)}=\psi_{\infty, {\rm GFL, q}}(\kappa; \hat{\bf r}_i,
{\bf d}_j)-\psi_{\infty, {\rm FL, q}}(\kappa; \hat{\bf r}_i,
{\bf d}_j)
\]
for the quadratically nonlinear point scatterers and 
\[
 P_\kappa^{(i, j)}=\psi_{\infty, {\rm GFL, c}}(\kappa; \hat{\bf r}_i,
{\bf d}_j)-\psi_{\infty, {\rm FL, c}}(\kappa; \hat{\bf r}_i,
{\bf d}_j)
\]
for the cubically nonlinear point scatterers. The purpose of taking the
difference is to avoid the possibility that the far-field pattern from
the point scatterers may dominate the far-field pattern from the extended
scatterers.

\section{Numerical experiments}

In this section, we discuss the implementation of the direct scattering problem
and present some numerical experiments for the inverse scattering problem. In
all of the following examples, the extended scatterer is a five-leaf shaped
obstacle and can be parameterized, up to a shift and rotation, by 
\begin{equation}\label{eo}
 {\bf r}(t)=r(t)(\cos t, \sin t), \quad r(t)=2+0.5\cos(5t), 
\end{equation}
where $t\in[0, 2\pi]$ is the parameter. For convenience, we summarize some of
the parameters used in the numerical experiments in Table \ref{tab-1}. The
resulting system of equations are obtained after discretizing the boundary of
the extended scatterers. The total number of unknowns is $N_{\rm point}+N_{\rm
boundary}$ for linear point scatterers and $2(N_{\rm point}+N_{\rm boundary})$
for nonlinear point scatterers.

\begin{table}
\caption{Parameters used in the numerical experiments.}
\begin{center}
\begin{tabular}{ l | l }
\hline
\hline
$N_{\rm point}$ & number of point scatterers  \\
\hline
$N_{\rm boundary}$ & number of points to discretize the boundary of the
extended scatterer(s)\\
\hline
$N_{\rm direction}$ & number of incident and observation directions\\
\hline
$N_{\rm sampling}$ & number of sampling points along the $x$- and
$y$-direction\\
\hline
$T_{\rm invert}$ & time (in seconds) to invert (factorize) the scattering
matrix\\
\hline
$T_{\rm solver}$ & time (in seconds) to solve the linear system for one
incidence\\
\hline
$T_{\rm ffp}$ & time (in seconds) to evaluate the far-field patterns\\
\hline
$T_{\rm NUFFT}$ & time (in seconds) to apply the NUFFT to evaluate the imaging
function\\
\hline
\hline
\end{tabular}
\end{center}
\label{tab-1}
\end{table}

\subsection{Direct scattering problem solver}

Let $\boldsymbol\phi=(\phi_1, \dots, \phi_m)^\top, \boldsymbol\phi_{\rm
inc}=(\phi_{\rm inc}({\bf r}_1), \dots, \phi_{\rm inc}({\bf r}_m))^\top$.
Define an $m\times m$ matrix
\[
 \mathcal{A}_\kappa=\begin{bmatrix}
                     1 & -\sigma_2 G_\kappa({\bf r}_1, {\bf r}_2) &\cdots& 
-\sigma_m G_\kappa({\bf r}_1, {\bf r}_m)\\
-\sigma_1 G_\kappa({\bf r}_2, {\bf r}_1) & 1 & \cdots & -\sigma_m G_\kappa({\bf
r}_2, {\bf r}_m)\\
\vdots & \vdots & \ddots & \vdots\\
-\sigma_1 G_\kappa({\bf r}_m, {\bf r}_1)&-\sigma_2 G_\kappa({\bf r}_m, {\bf
r}_2) & \cdots & 1  
                    \end{bmatrix}
\]
and three linear operators
\begin{align*}
& \mathcal{M}_\kappa u=\left(\int_\Gamma G_\kappa({\bf r}_1, {\bf
r}')u({\bf r}'){\rm d}s({\bf r}'), \dots,  \int_\Gamma G_\kappa({\bf r}_m, {\bf
r}')u({\bf r}'){\rm d}s({\bf r}')\right)^\top,\\
& (\mathcal{N}_\kappa \boldsymbol{\phi})({\bf r})=-\sum_{k=1}^m \sigma_k\phi_k
(\partial_\nu-{\rm i}\eta) G_\kappa({\bf r}, {\bf r}_k),\\
&(\mathcal{K}_\kappa u)({\bf r})=\frac{1}{2}u({\bf r})+\int_\Gamma
(\partial_\nu-{\rm i}\eta)G_\kappa({\bf r}, {\bf r}')u({\bf r}'){\rm d}s({\bf
r}').
\end{align*}
The generalized Foldy--Lax formulation of linear point scatterers
\eqref{gfl-lps} and \eqref{gfl-li} can be written as the operator form:
\begin{equation}\label{ls-l}
\begin{bmatrix} 
\mathcal{A}_\kappa & \mathcal{M}_\kappa \\ 
\mathcal{N}_\kappa & \mathcal{K}_\kappa 
\end{bmatrix}
\begin{bmatrix} 
\boldsymbol\phi  \\ 
\varphi  
\end{bmatrix} = 
\begin{bmatrix} 
\boldsymbol\phi_{\rm inc}  \\ 
\varphi_{\rm inc}
\end{bmatrix},
\end{equation}
where $\varphi=\partial_\nu\phi({\bf r})$ and $\varphi_{\rm
inc}=(\partial_\nu-{\rm i}\eta)\phi_{\rm inc}({\bf r})$. In the
discretization, the surface of the extended obstacle $\Gamma$ is discretized
by a set of uniform points in the parameter space $t$; the singular integral
is evaluated by the Alpert quadrature \cite{a-sisc99}; the boundary integral
equations are solved by the Nystr\"om method.

There are three approaches to solve the linear system \eqref{ls-l}:

\begin{enumerate}

\item Direct solver. Apply the LU factorization from the Lapack library
and parallelize it by OpenMP on a multicore workstation;

\item Iterative solver. Apply GMRES to the whole system and accelerate
the matrix vector product by the fast multipole method (FMM)\cite{GR1987}.

\item Hybrid method. Assume that the number of points to discretize the
extended obstacle is relatively small compared to the number of point
scatterers. First is to invert $\mathcal{K}_\kappa$ directly and then solve the
Schur complement of \eqref{ls-l} by an iterative method, i.e., solve
iteratively with the FMM acceleration of the linear system
\[
(\mathcal{A}_\kappa-\mathcal{M}_\kappa\mathcal{K}_\kappa^{-1}\mathcal{N}
_\kappa)\boldsymbol\phi=\boldsymbol\phi_{\rm inc}-\mathcal{M}_\kappa
\mathcal{K}_\kappa^{-1}\varphi_{\rm inc}.
\]
\end{enumerate}              

To investigate the three different methods, we solve the scattering problem for
the extended obstacle in \eqref{eo}, which is surrounded by a group of linear point
scatterers. Table \ref{tab-2} shows the numerical performance for the
three different methods. Obviously, the direct solver is the best
option for the linear problem. It solves the system very rapidly with
parallelization, and it is independent of the location of the point scatterers
and the wavenumber $\kappa$. Both of the iterative methods fail if the point
scatterers are randomly and densely distributed in a specific area. In addition,
to solve the inverse problem, the direct problem \eqref{ls-l} has to be
solved many times with different right hand sides in order to construct the
response matrix, which corresponds to the far-field patterns for different
incident directions. The advantage of the direct solver is clear: we only need
to invert \eqref{ls-l} once and apply matrix vector product for different
right hand sides, which can also be done by parallelization. 

\begin{table}
\caption{Time (in seconds) to solve the linear system \eqref{ls-l}
on an HP workstation.}
\begin{center}
\begin{tabular}{ c c c c c }
\hline
\hline
$N_{\rm point}$ & $N_{\rm boundary}$ & Method 1 & Method 2 & Method 3   \\
\hline
$1000$ & $600$ & $0.16$ & $1.42$  & $0.89$ \\
\hline
$10000$ & $600$ & $8.9$ & fail to converge  & fail to converge
\\
\hline
\hline
\end{tabular}
\end{center}
\label{tab-2}
\end{table}

Next we discuss the solver for the generalized Foldy--Lax of nonlinear point
scatterers. We only describe the steps for the quadratically nonlinear point
scatterers since they are similar to the cubically nonlinear point scatterers. 

Let $\boldsymbol\phi^{(j)}=(\phi_1^{(j)}, \dots, \phi_m^{(j)})^\top$ and
$\varphi^{(j)}=\partial_\nu\phi^{(j)}({\bf r})$ be the external field acting on
the point scatterers and the normal derivative of the total field on the
boundary of the extended obstacle at the wavenumber $\kappa_j$, respectively.
The generalized Foldy--Lax for nonlinear point scatterers can be written as
\begin{equation}\label{ls-n}
\setlength{\dashlinegap}{2pt}
\left[\begin{array}{c c : c c} 
\mathcal{D}_{1 1}& \mathcal{D}_{1 2} & \mathcal{M}_{\kappa_1} & 0 \\
\mathcal{D}_{2 1} & \mathcal{D}_{2 2}& 0 & \mathcal{M}_{\kappa_2}\\
\hdashline
\mathcal{H}_{1 1} & \mathcal{H}_{1 2} & \mathcal{K}_{\kappa_1} & 0\\
\mathcal{H}_{2 1} & \mathcal{H}_{2 2} & 0 &\mathcal{K}_{\kappa_2}
\end{array}\right]
\left[\begin{array}{c} 
\boldsymbol\phi^{(1)} \\
\boldsymbol\phi^{(2)}\\
\hdashline
\varphi^{(1)}\\
\varphi^{(2)}
\end{array}\right]= 
\left[
\begin{array}{c}
\boldsymbol\phi_{\rm inc}  \\ 
0\\
\hdashline
\varphi_{\rm inc}\\
0
\end{array}
\right],
\end{equation}
where $\mathcal{D}_{i j}$ represents the nonlinear interaction between
$\boldsymbol\phi^{(1)}$ and $\boldsymbol\phi^{(2)}$ at the point scatterers,
$\mathcal{H}_{i j}$ is the nonlinear interaction from the point scatterers to
the extended obstacle, $\mathcal{M}_{\kappa_j}$ denotes the linear interaction
from the extended obstacle to the point scatterers at the wavenumber
$\kappa_j$, and $\mathcal{K}_{\kappa_j}$ is the linear interaction for the
extended obstacle.

Due to the large number of unknowns and nonlinearity, neither the direct solver
nor the iterative method is applicable to the nonlinear system \eqref{ls-n}.
Assuming that the number of point scatterers is relatively small, we propose an
efficient nonlinear solver, which can be be applied to the Schur complement of
\eqref{ls-n}. Specifically, we invert $\mathcal{K}_{\kappa_1}$ and
$\mathcal{K}_{\kappa_2}$ directly, and then solve the following nonlinear
system:
\begin{align}\label{sc-n}
\left\{
\begin{bmatrix} 
\mathcal{D}_{1 1} & \mathcal{D}_{1 2} \\ 
\mathcal{D}_{2 1} & \mathcal{D}_{2 2} \\ 
\end{bmatrix}
-
\begin{bmatrix} 
\mathcal{M}_{\kappa_1}\mathcal{K}^{-1}_{\kappa_1} & 0 \\ 
0 & \mathcal{M}_{\kappa_2}\mathcal{K}^{-1}_{\kappa_2} \\ 
\end{bmatrix}
\begin{bmatrix} 
\mathcal{H}_{1 1} & \mathcal{H}_{1 2} \\ 
\mathcal{H}_{2 1} & \mathcal{H}_{2 2} \\ 
\end{bmatrix}
\right\}
\begin{bmatrix} 
\boldsymbol\phi^{(1)}  \\ 
\boldsymbol\phi^{(2)} \\
\end{bmatrix}\notag\\
=\begin{bmatrix} 
\boldsymbol\phi_{\rm inc}  \\ 
0  \\
\end{bmatrix}
-
\begin{bmatrix} 
\mathcal{M}_{\kappa_1}\mathcal{K}^{-1}_{\kappa_1} & 0 \\ 
0 & \mathcal{M}_{\kappa_2}\mathcal{K}^{-1}_{\kappa_2} \\ 
\end{bmatrix}
\begin{bmatrix}
\varphi_{\rm inc} \\
0
\end{bmatrix}.
\end{align}
Finally, we apply a trust region Newton type method to solve the
resulting nonlinear system \eqref{sc-n}.


\subsection{Linear point scatterers}

We show the imaging results of extended scatterers which are surrounded by a
group of linear point scatterers. These point scatterers are randomly and
uniformly distributed in the annulus $\{{\bf r}\in\mathbb R^2: 10 <
|{\bf r}|<11\}$. The scattering coefficients are $\sigma_k=0.5$ for all the
point scatterers. The numerical performance of the following two examples is
shown in Table \ref{tab-3}. It is clear to note the the proposed method is not
only efficient for the direct problem simulation but also for the inverse
problem imaging. 

\begin{table}
\caption{Results for imaging two extended scatterers surrounded by linear point
scatterers.} 
\begin{center}
\begin{tabular}{ c c c c c c c c c c c }
\hline
\hline
& $\kappa$ &$N_{\rm point}$ & $N_{\rm boundary}$ & $N_{\rm direction}$ &
$N_{\rm sampling}$ & $T_{\rm invert}$ & $T_{\rm sampling}$ & $T_{\rm ffp}$ &
$T_{\rm NUFFT}$ \\
\hline
Example 1 & 10 & 1000 & 600 & 360 & 500 & 7.95e-2 &
2.56e-3 & 2.23e-2 & 2.46e-1 \\
\hline
Example 2 & 50 & 1000 & 4800 & 1800 &  500 & 1.61 &
2.17e-2  & 3.49e-1 & 3.70\\
\hline
\hline
\end{tabular}
\end{center}
\label{tab-3}
\end{table}
 
\subsubsection{Example 1} 

This example is to image two extended scatterers which are surrounded by
1000 linear point scatterers at the wavenumber $\kappa=10$. The imaging result
is shown in Figure \ref{fig-2}(a). The extended scatterers are well
reconstructed except for the parts where these two scatterers are close. 

\subsubsection{Example 2} 

\begin{figure}
\centering
\includegraphics[width=0.4\textwidth]{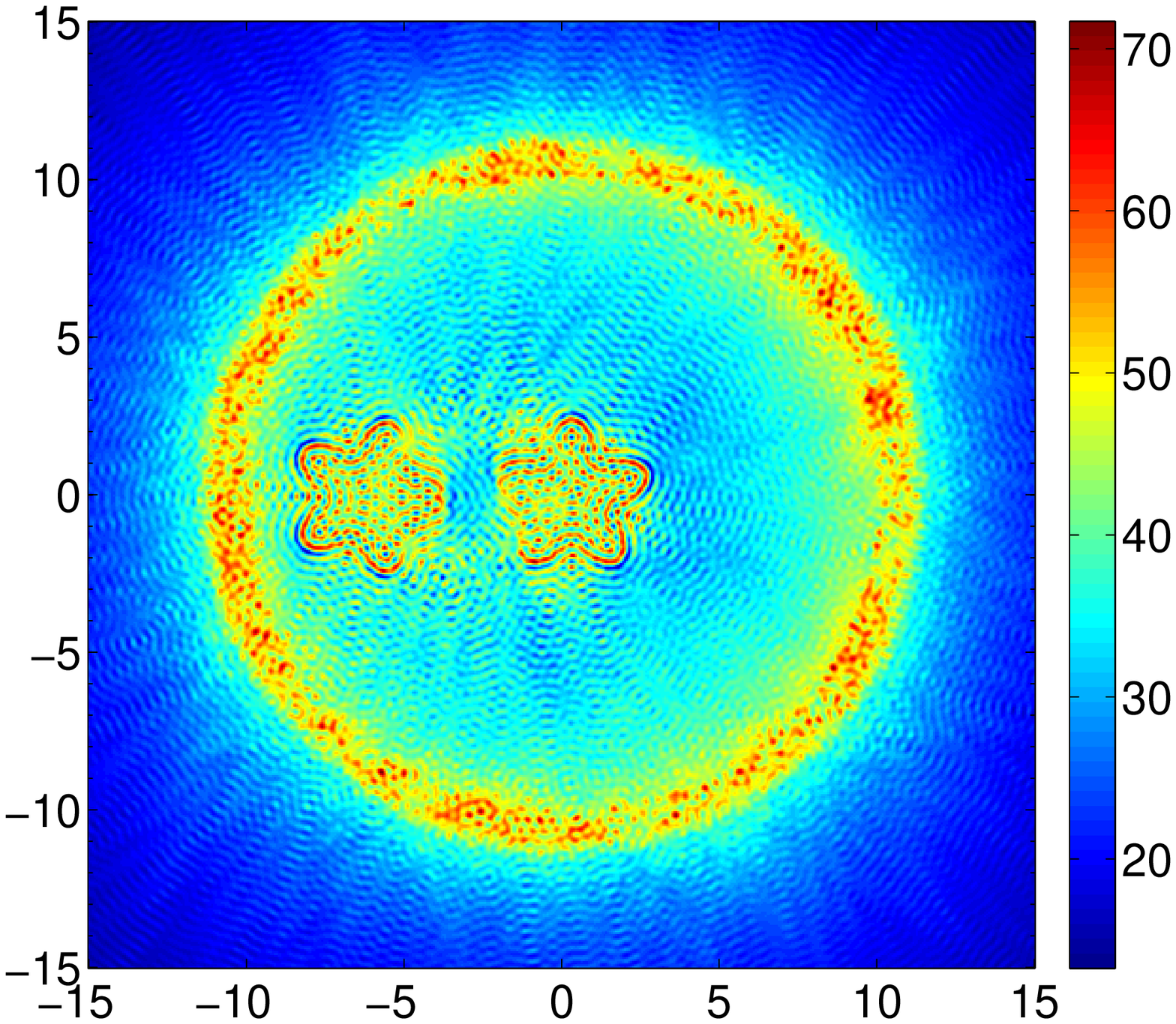}
\includegraphics[width=0.4\textwidth]{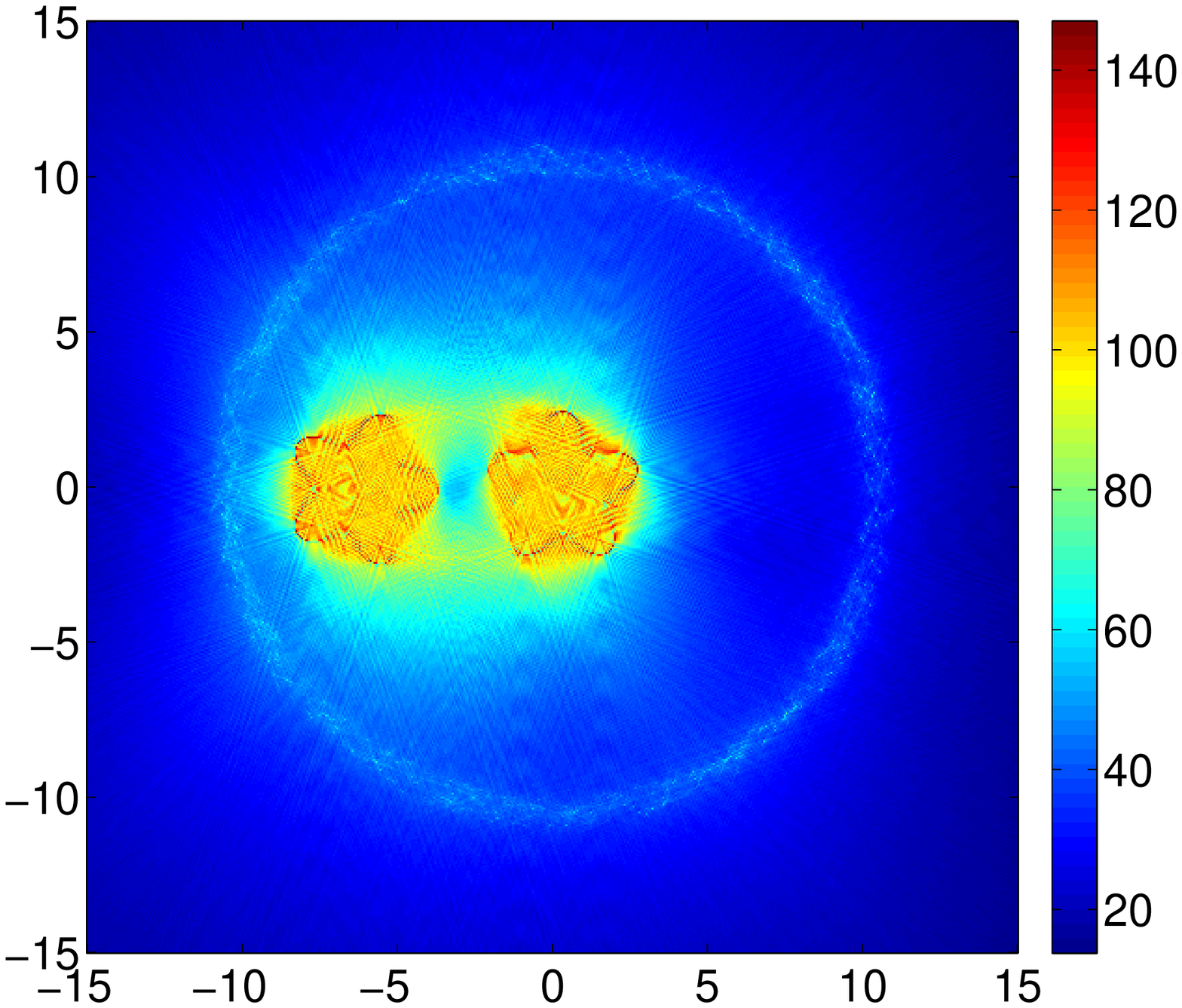}\\
(a)\hspace{6cm} (b)
\caption{Imaging of two extended scatterers surrounded by
1000 linear point scatterers. (a) Example 1: $\kappa=10$;
(b) Example 2: $\kappa=50$.}
\label{fig-2}
\end{figure}

To show the influence of the wavenumber on the imaging resolution, we take the
same extended and point scatterers as those in Example 1, but we use the
wavenumber $\kappa=50$. The imaging result is shown in Figure \ref{fig-2}(b),
which has a better resolution than Figure \ref{fig-2}(a) does. These two
scatterers are well reconstructed even for the parts where they are close. As
is expected, higher wavenumber can capture finer structures.

\subsection{Quadratically nonlinear point scatterers}

We consider the imaging of extended scatterers with two quadratically nonlinear
point scatterers. The nonlinear scattering coefficients are
$\sigma_{k,1}^{(1)}=\sigma_{k,2}^{(1)}=0.5$
and $\sigma_{k,1}^{(2)}=\sigma_{k,2}^{(2)}=0.4$. The numerical performance is
shown in Table \ref{tab-4} for the following two examples. As can be seen, the
proposed method is also efficient for the nonlinear direct problem simulation. 

\begin{table}
\caption{Results for imaging the extended scatterers surrounded by quadratically
nonlinear point scatterers.}
\begin{center}
\begin{tabular}{ c c c c c c c c c c }
\hline
\hline
& $\kappa$ &$N_{\rm point}$ & $N_{\rm boundary}$ & $N_{\rm direction}$ &
$N_{\rm sampling}$ & $T_{\rm invert}$ & $T_{\rm solver}$ & $T_{\rm ffp}$ &
$T_{\rm NUFFT}$ \\
\hline
Example 3 & 2 & 2 & 600 & 360 & 500 & 1.22e-3 & 8.52e-3  &
1.39e-2 & 3.39e-1 \\
\hline
Example 4 & 5 & 2 & 1200 & 360 & 500 & 1.90e-1 & 3.99e-3 
& 1.96e-2 & 3.84e-1\\
\hline
\hline
\end{tabular}
\end{center}
 \label{tab-4}
\end{table}

\subsubsection{Example 3}

\begin{figure}
\centering
\includegraphics[width=0.4\textwidth]{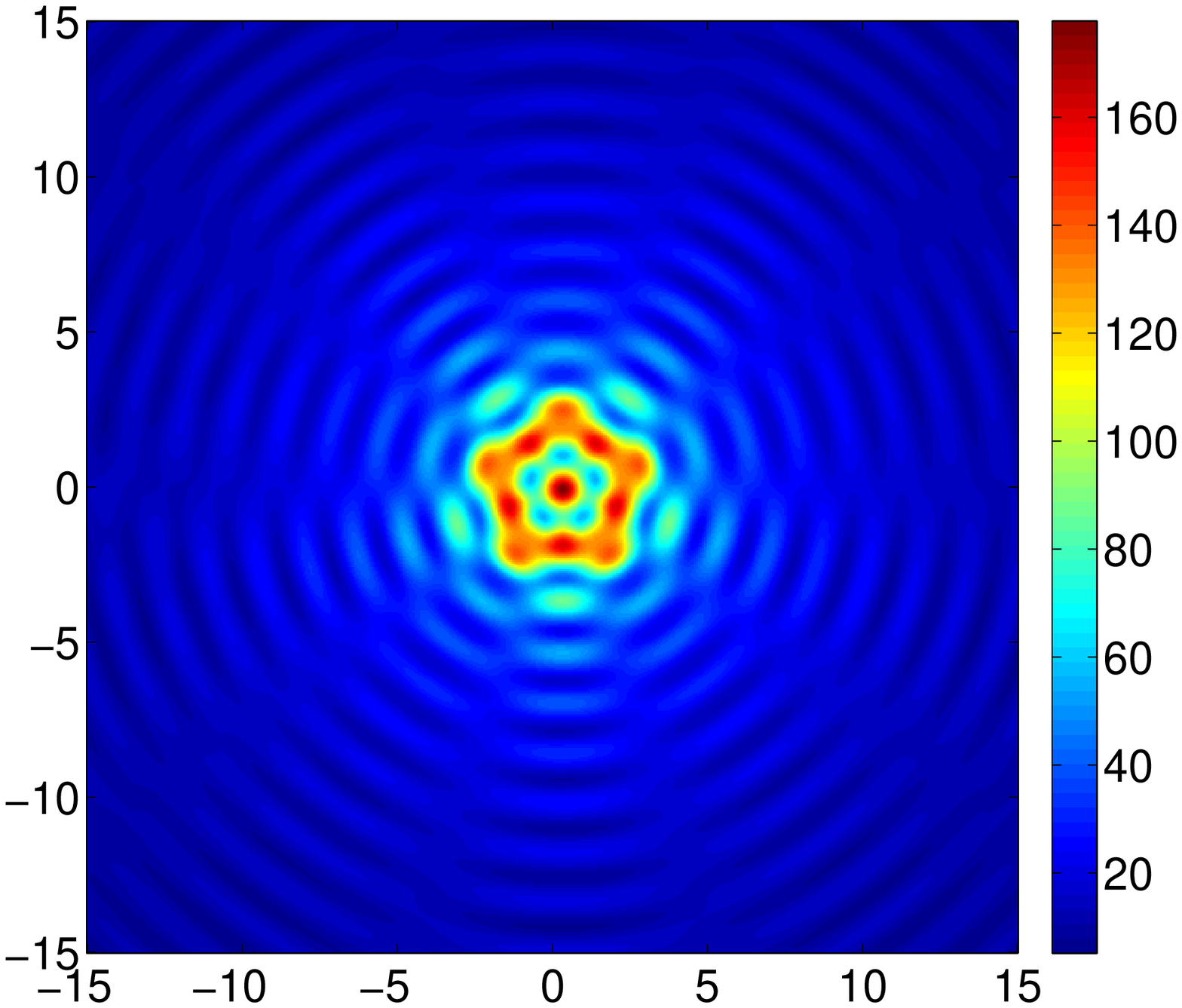}
\includegraphics[width=0.4\textwidth]{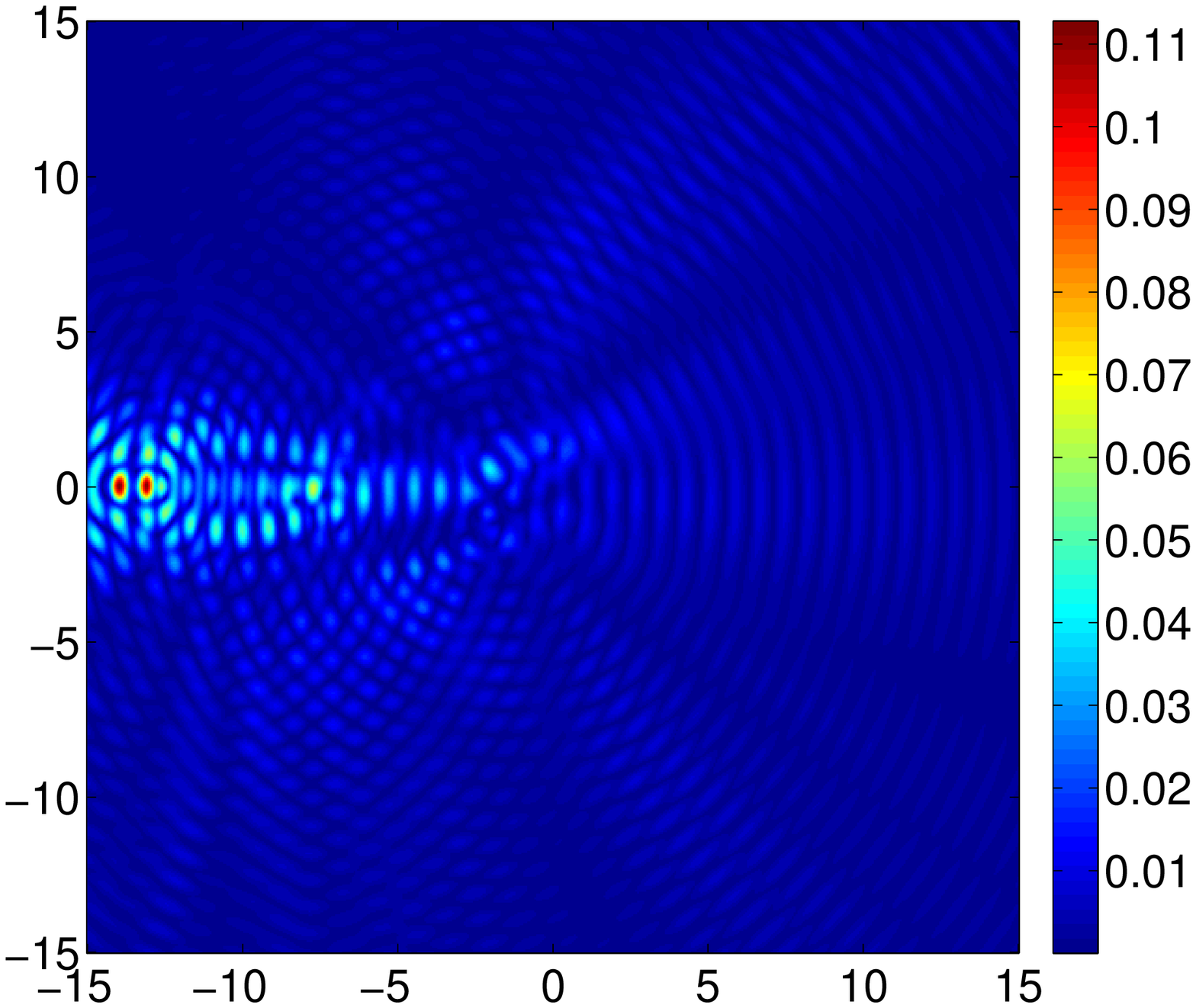}\\
(a)  \hspace{6cm} (b) 
\caption{Example 3: Imaging of one extended scatterer with
two fixed quadratically nonlinear point scatterers. (a)
Imaging with $\kappa_1=2$; (b) Imaging with $\kappa_2=4$.}
\label{fig-3}
\end{figure}

\begin{figure}
\centering
\includegraphics[width=0.4\textwidth]{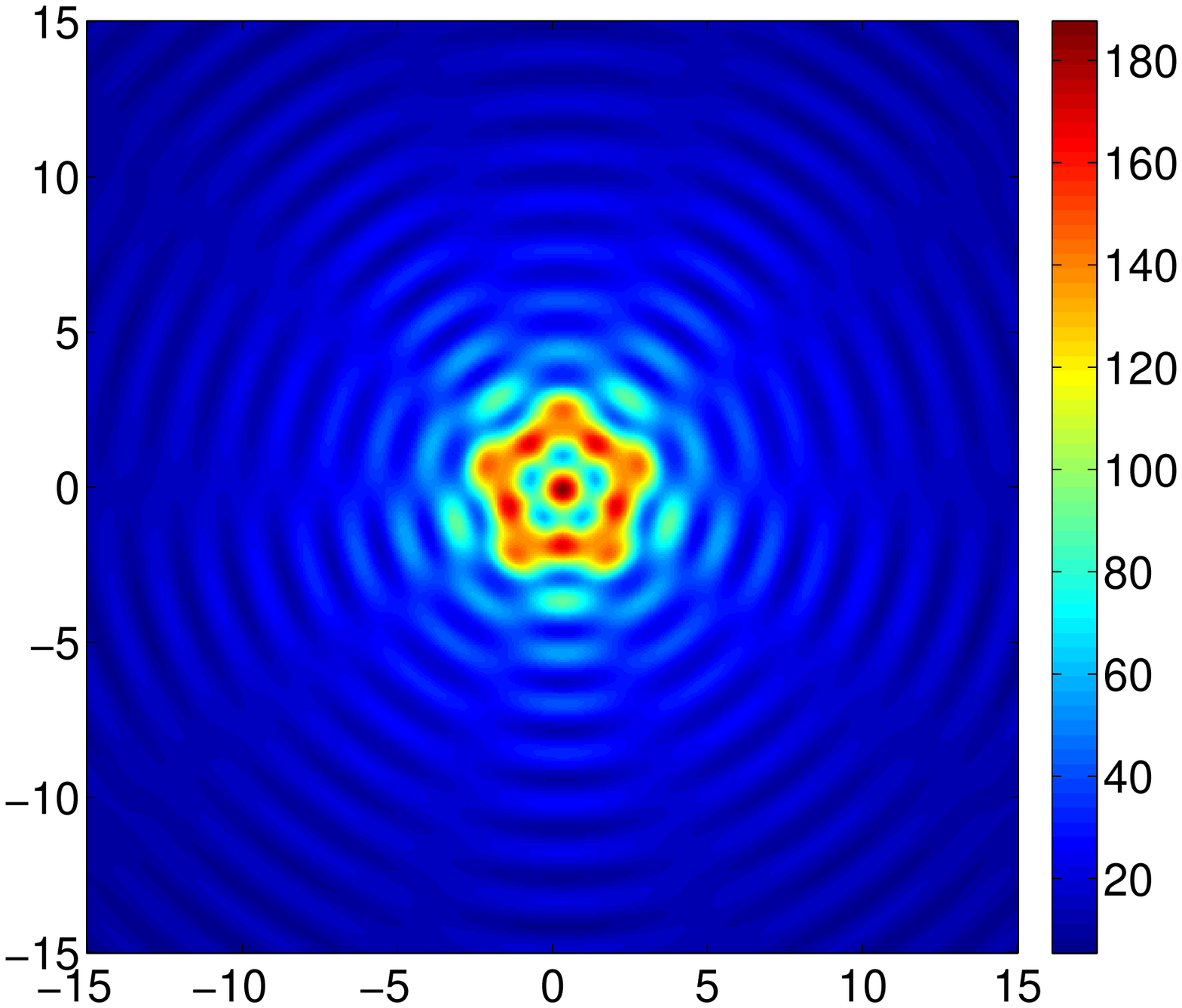}
\includegraphics[width=0.4\textwidth]{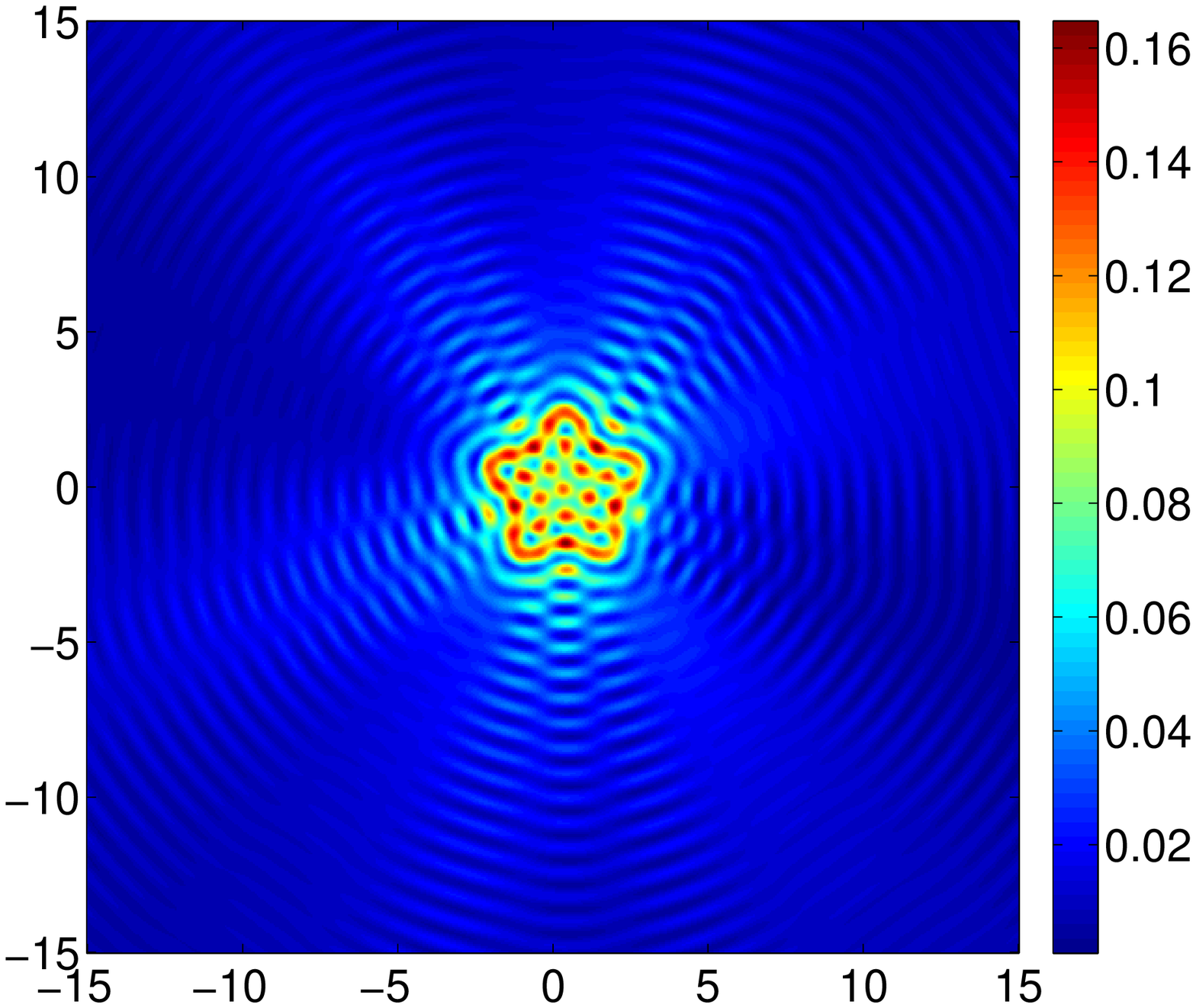}\\
(a)  \hspace{6cm} (b) 
\caption{Example 3: Imaging of one extended scatterer with
two moving quadratically nonlinear point scatterers. (a)
Imaging with $\kappa_1=2$; (b) Imaging with $\kappa_2=4$.}
\label{fig-4}
\end{figure}

This example is to image one extended scatterers with two quadratically
nonlinear point scatterers. The wavenumber of the incident wave is
$\kappa=2$. First we consider the case when the two point scatterers are fixed
at the location $(-13, 0)$ and $(-14, 0)$, respectively. The imaging result is
shown in Figure \ref{fig-3}. The linear wave can reconstruct the extended
scatterer but with poor resolution. The wave from the second harmonic
generation cannot reconstruct the extended scatterer. As is described in
section 4.4, we change the location of the two point scatterers on the circle
with radius 13 and 14 aligned with the incident direction. For instance, if the
angle of the incidence is $\theta=\frac{\pi}{3}$, the two point scatterers are
located at $(13\cos\theta, 13\sin\theta)$ and $(14\cos\theta, 14\sin\theta)$.
The imaging result is shown in Figure \ref{fig-4}. The linear wave yields
almost the same imaging result as the fixed point scatterers does. But the
wave from the second harmonic generation gives a much better imaging result due
to the doubled frequency. 

\subsubsection{Example 4}

\begin{figure}
\centering
\includegraphics[width=0.4\textwidth]{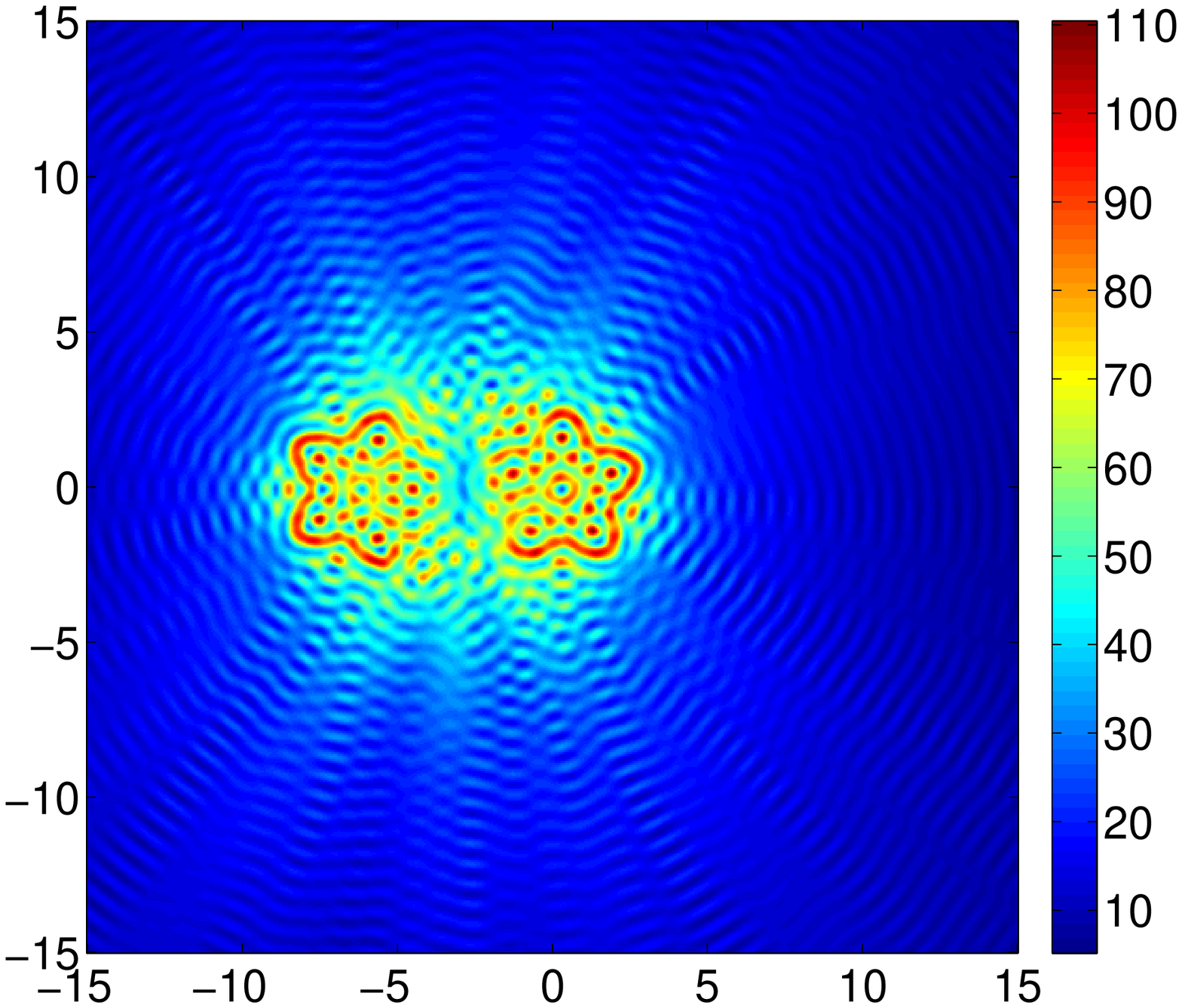}
\includegraphics[width=0.4\textwidth]{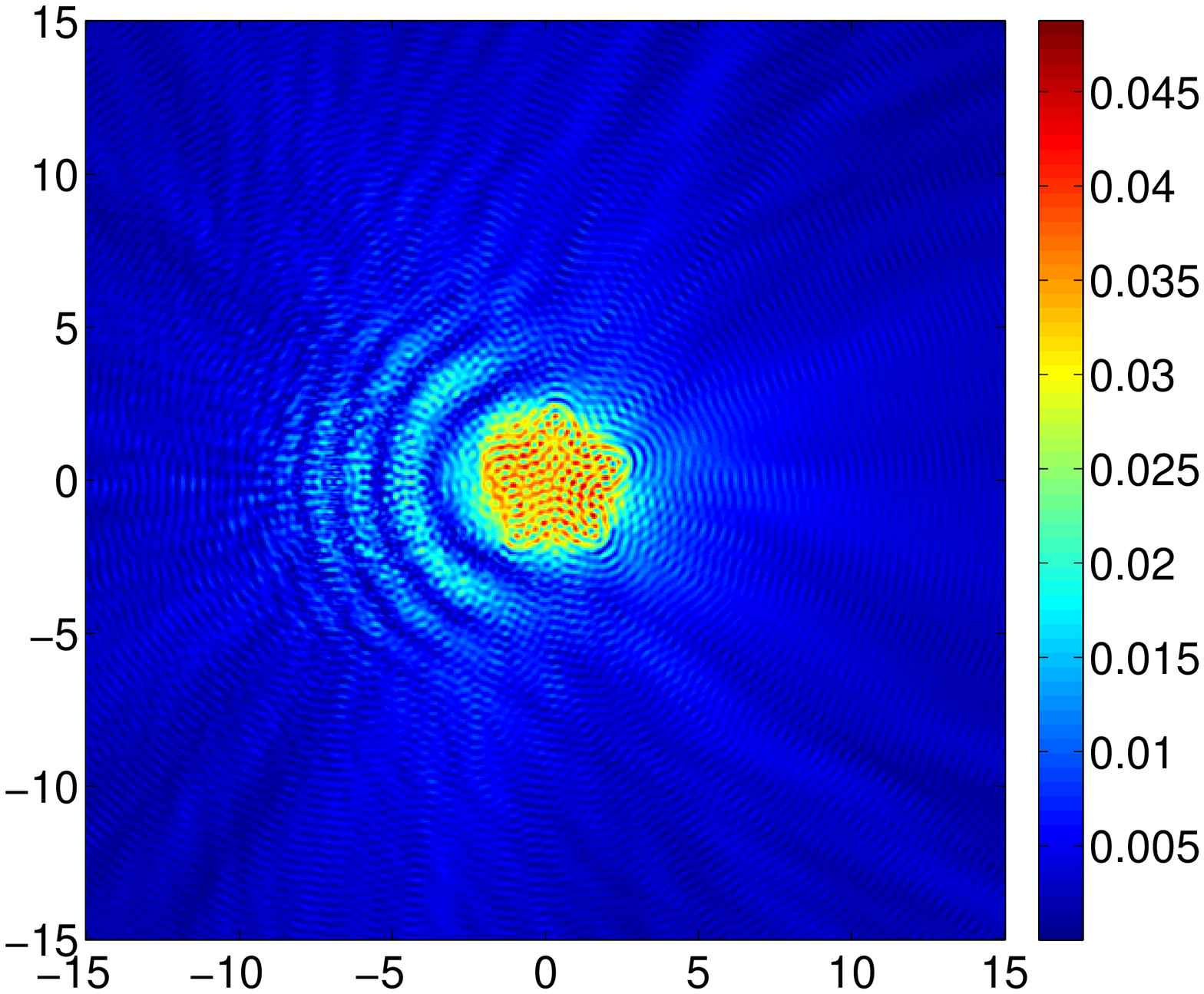}\\
(a)  \hspace{6cm} (b) 
\caption{Example 4: Imaging of two extended scatterers with two
quadratically nonlinear point scatterers close by. (a) Imaging with $\kappa_1=5$; (b)
Imaging with $\kappa_2=10$.}
\label{fig-5}
\end{figure}

\begin{figure}
\centering
\includegraphics[width=0.4\textwidth]{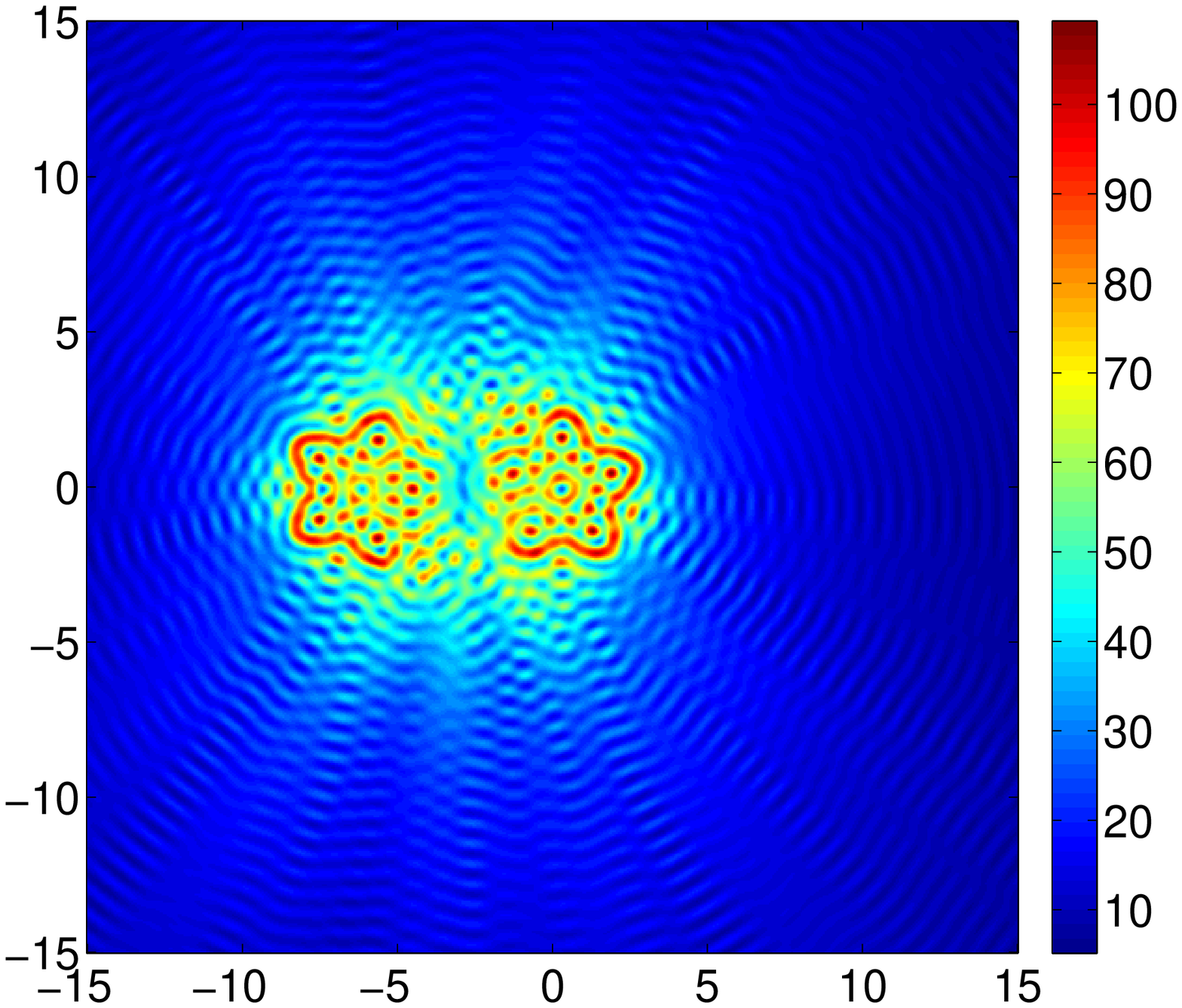}
\includegraphics[width=0.4\textwidth]{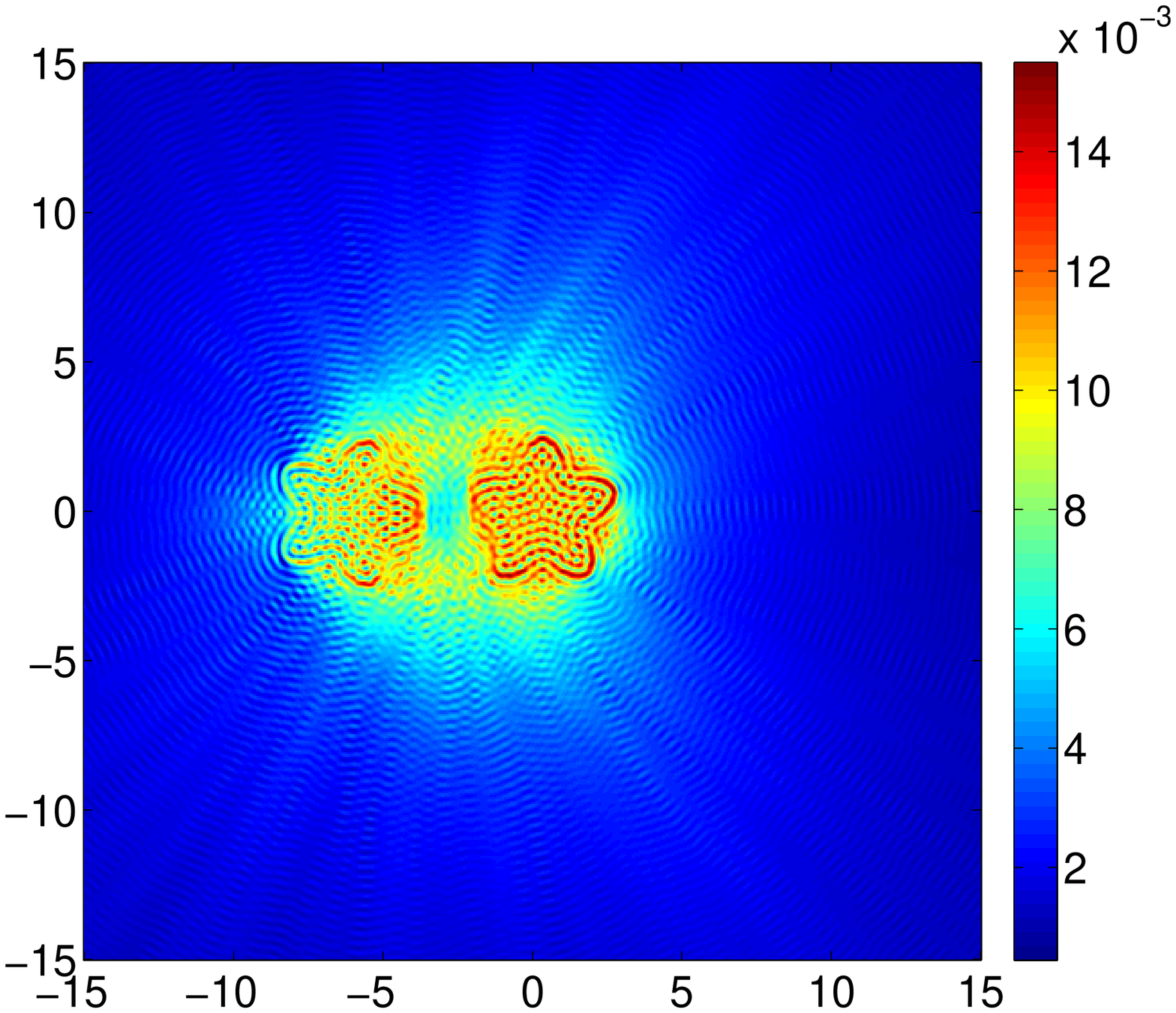}\\
(a)  \hspace{6cm} (b) 
\caption{Example 4: Imaging of two extended scatterers with two
quadratically nonlinear point scatterers far away. (a) Imaging with $\kappa_1=5$; (b)
Imaging with $\kappa_2=10$.}
\label{fig-6}
\end{figure}

This example is to image two extended scatterers with two point scatterers
located on the circles of radii $|{\bf r}|=13$ and $|{\bf r}|=14$. The
wavenumber of the incident wave is $\kappa=5$. The imaging result is shown in
Figure \ref{fig-5}. The linear wave can still produces a reasonable imaging
result. But the nonlinear wave fails to identify the two extended scatterers. 
Next we keep everything else the same except moving away the two point
scatterers to the circles of radii $|{\bf r}|=130$ and $|{\bf r}|=131$. The
imaging result is shown in Figure \ref{fig-6}. It can be seen that the nonlinear
wave can clearly identify the two extended scatterers. The reason for this is
clear: when the point scatterers are far away, their generated wave interacted
with the extended obstacles is more like a plane wave incidence, which can be
better resolved by the illumination vectors.   

\subsection{Cubically nonlinear point scatterers}

Finally, we investigate the performance of using cubically nonlinear point
scatterers. We put two cubically nonlinear point scatterers around the extended
scatterers. The scattering coefficients
are $\sigma_{k,1}^{(1)}=\sigma_{k,2}^{(1)}=0.5$ and
$\sigma_{k,1}^{(3)}=\sigma_{k,2}^{(3)}=\sigma_{k,3}^{(3)}=0.4$. The numerical
performance is summarized in Table \ref{tab-5}. 

\begin{table}
\caption{Results for imaging the extended scatterers surrounded by
cubically nonlinear point scatterers.} 
\centering
\begin{tabular}{ c c c c c c c c c c }
\hline
\hline
& $\kappa$ &$N_{\rm point}$ & $N_{\rm boundary}$ & $N_{\rm direction}$ &
$N_{\rm sampling}$ & $T_{\rm invert}$ & $T_{\rm solver}$ & $T_{\rm ffp}$ &
$T_{\rm NUFFT}$ \\
\hline
Example 5 & 2 & 2 & 600 & 360 & 500 & 1.24e-3 & 1.00e-2  &
1.22e-2 & 4.33e-1 \\
\hline
Example 6 & 5 & 2 & 1200 & 360 & 500 & 4.22e-3 & 1.91e-2 
& 2.11e-2 & 4.41e-1\\
\hline
\hline
\end{tabular}
\label{tab-5}
\end{table}

\subsubsection{Example 5}

\begin{figure}
\centering
\includegraphics[width=0.4\textwidth]{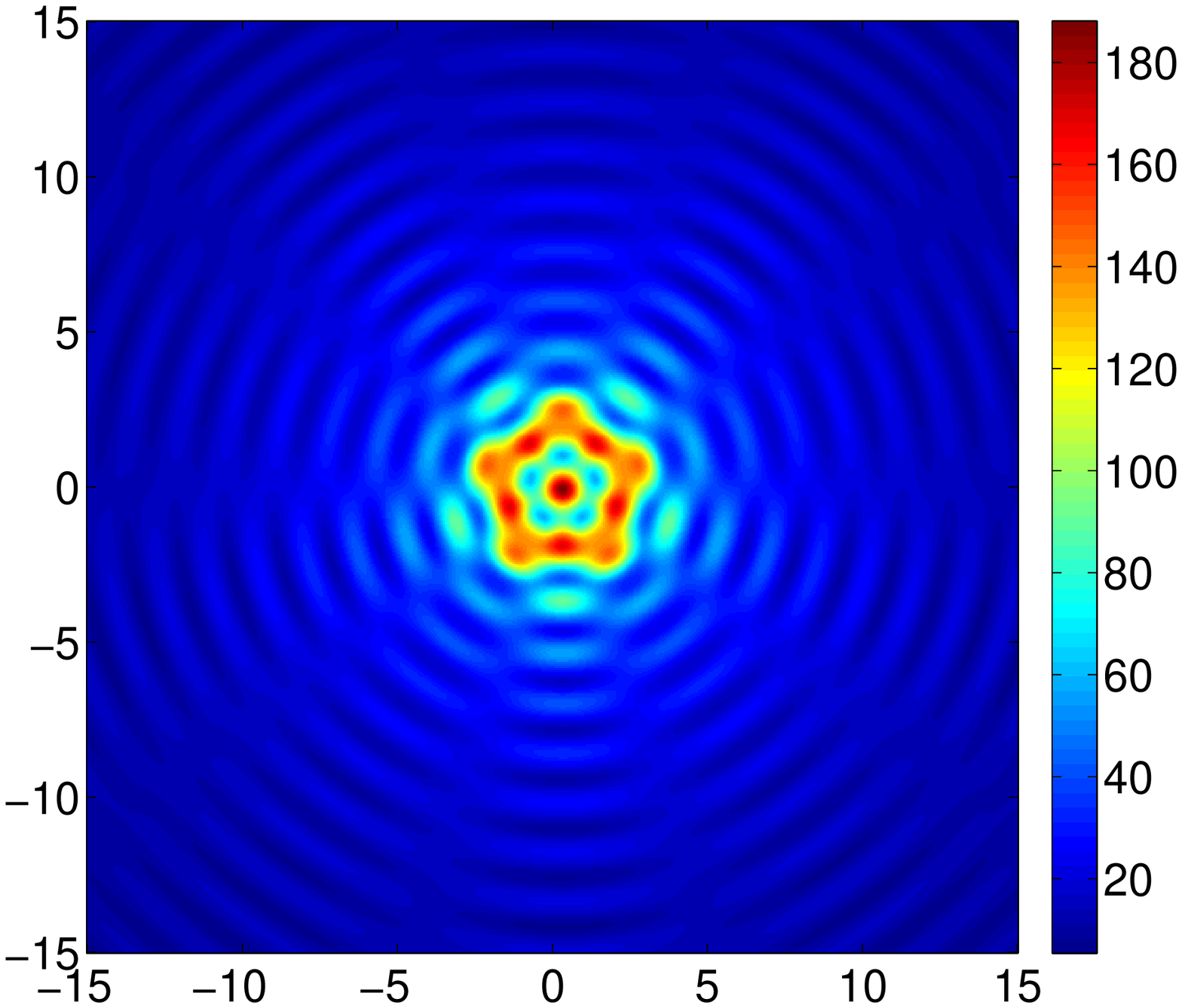}
\includegraphics[width=0.4\textwidth]{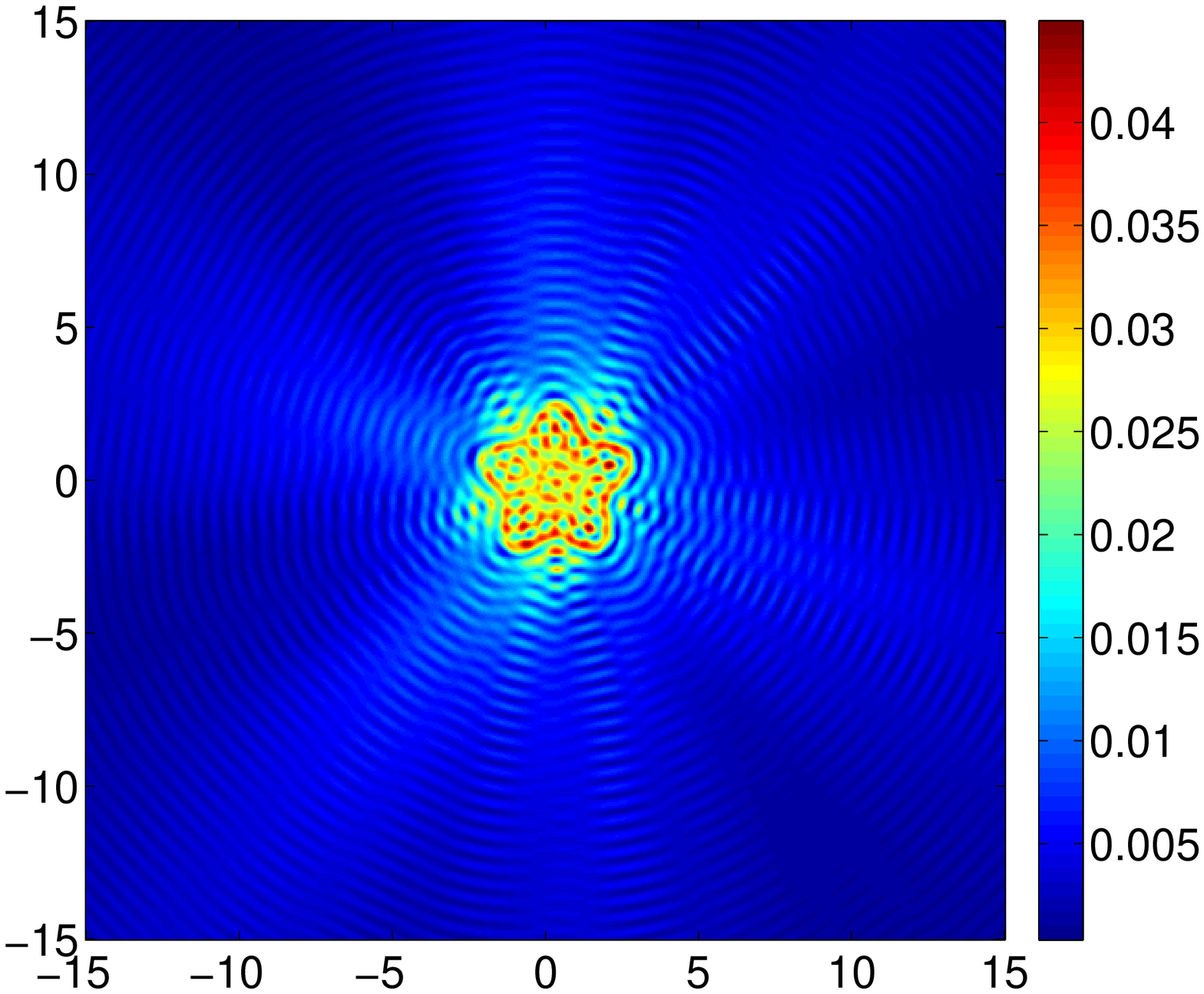}\\
(a)  \hspace{6cm} (b) 
\caption{Example 5: Imaging of two extended scatterers with two
quadratically nonlinear point scatterers. (a) Imaging with $\kappa_1=2$; (b)
Imaging with $\kappa_3=6$.}
\label{fig-7}
\end{figure}

This example is to image one extended scatterer with two point scatterers
located at the circles with radii $|{\bf r}|=13$ and $|{\bf r}|=14$,
respectively. The wavenumber of the incidence is $\kappa=2$. The imaging result
is shown in Figure \ref{fig-7}. Comparing Figure \eqref{fig-4}(b) with Figure
\eqref{fig-7}(b), we observe that the cubically nonlinear point scatterers can
produce a better resolution than the quadratically nonlinear point scatterers
does, which confirms once again that higher frequency produces finer
resolution. 

\subsubsection{Example 6}

\begin{figure}
\centering
\includegraphics[width=0.4\textwidth]{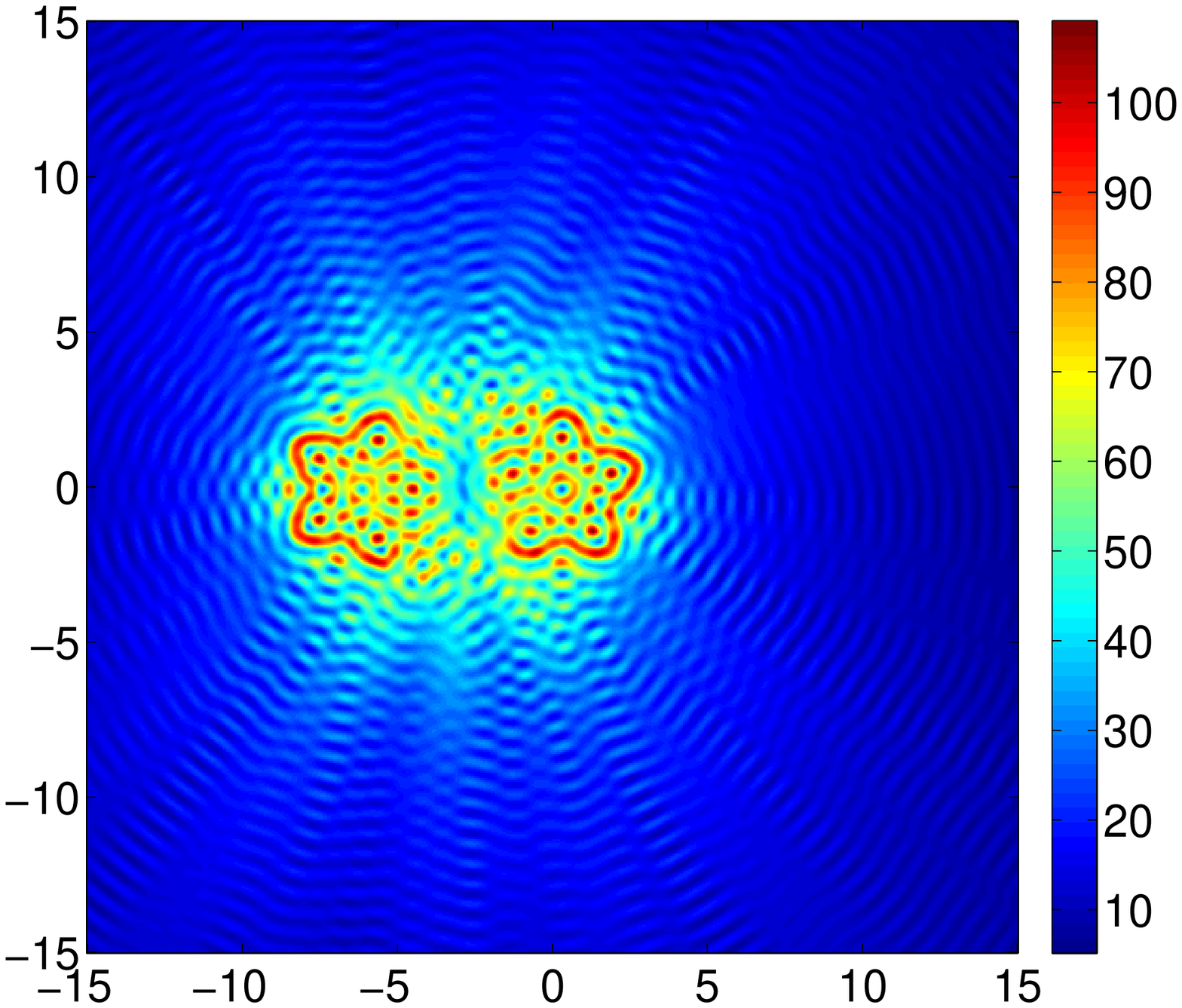}
\includegraphics[width=0.4\textwidth]{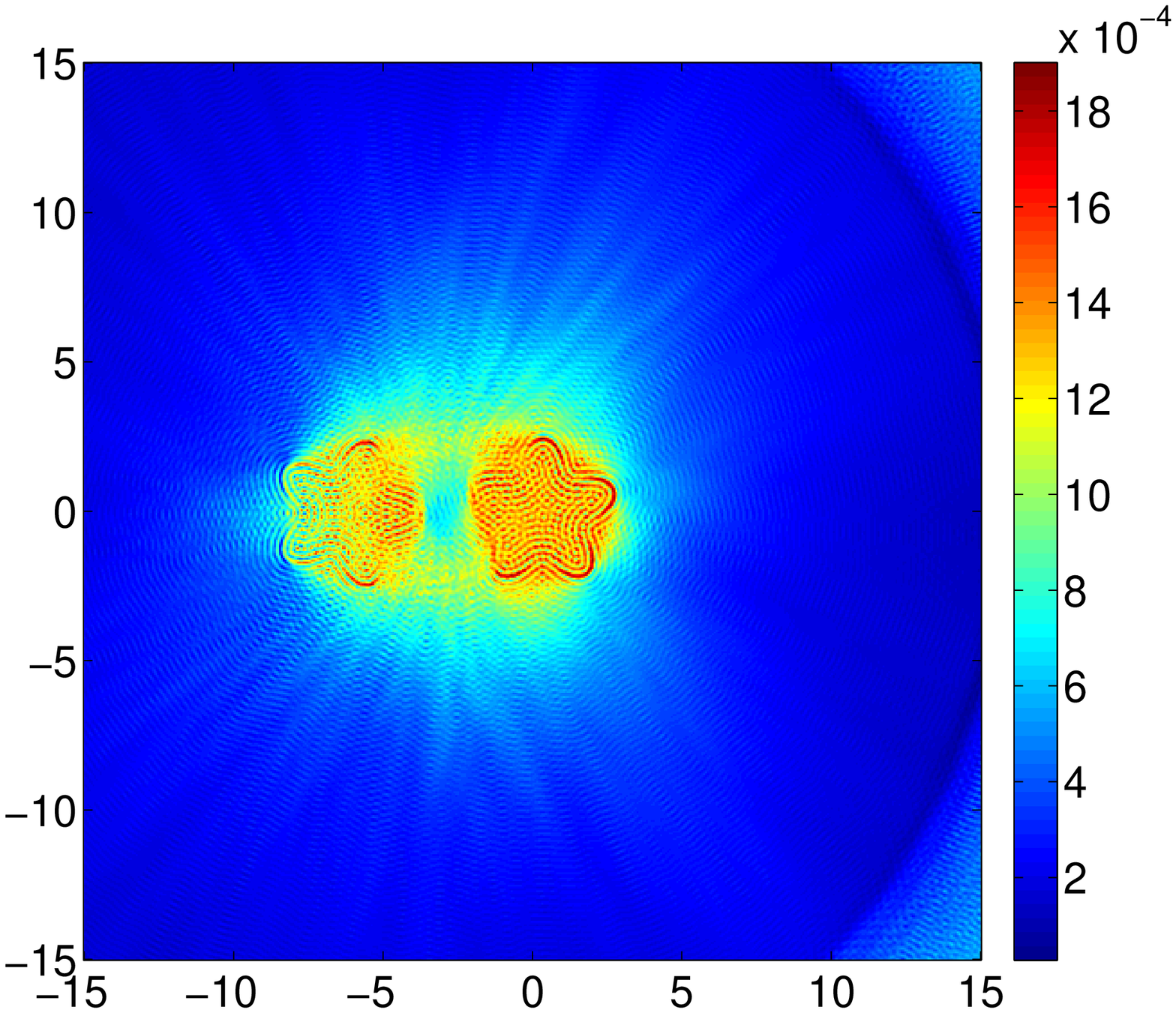}\\
(a)  \hspace{6cm} (b) 
\caption{Example 6: Imaging of two extended scatterers with two
quadratically nonlinear point scatterers. (a) Imaging with $\kappa_1=5$; (b)
Imaging with $\kappa_3=15$.}
\label{fig-8}
\end{figure}

This example is to image two extended scatterers with two point scatterers
located at the circles with radii $|{\bf r}|=13$ and $|{\bf r}|=14$,
respectively. The wavenumber of the incidence is $\kappa=5$. The imaging result
is shown in Figure \ref{fig-8}. Comparing Figure \eqref{fig-6}(b) with Figure
\eqref{fig-8}(b), we observe the same pattern that higher frequency wave
generates better resolved imaging result. 

\section{Conclusion}

We have presented a generalized Foldy--Lax formulation for the scattering by
the heterogeneous media consisting of linear or nonlinear point scatterers and
extended obstacles. A new imaging function is proposed and a fast direct imaging
method is developed for the inverse obstacle scattering problem. Using the
nonlinear point scatterers to excite high harmonic generation, enhanced imaging
resolution is achieved to reconstruct the extended obstacle. Our method shares
the attractive features of all the other direct imaging methods. In addition,
the evaluation is accelerated by using the FFT due to the special construction
of the imaging function. Numerical results show that the method is effective to
solve the inverse obstacle scattering problem. The proposed method can be
directly extended to solve the inverse obstacle scattering problem in three
dimensions. As an FFT-based method, a greater potential can show up in higher
dimensions. We are currently working on the three-dimensional problem and the
results will be reported elsewhere.  

\appendix

\section{Nonlinear wave equations}

In the scalar theory of electromagnetic fields, the scalar electric field
$u(\textbf{r},t)$ obeys the wave equation:
\begin{equation}\label{a-we}
\Delta u({\bf r},t)-\frac{1}{c^2}\partial^2_t u({\bf
r},t)=\frac{4\pi}{c^2}\partial^2_t P({\bf r},t),\quad{\bf r}\in\mathbb R^2,
\end{equation}
where $c>0$ is the light speed and $P({\bf r}, t)$ is the polarization density. 

We adopt the following Fourier and inverse Fourier transformation convention: 
\[
f(\textbf{r}, \omega)=\int_{\mathbb R} f(\textbf{r}, t)e^{{\rm i}\omega t}{\rm
d}t, \quad f(\textbf{r}, t) =\frac{1}{2\pi}\int_{\mathbb R} f(\textbf{r},
\omega)e^{-{\rm i}\omega t}{\rm d}\omega.
\]
Note that if $f(\textbf{r},t)$ is a real-valued function, then
$f(\textbf{r},-\omega)=\bar{f}(\textbf{r},\omega)$, where the bar is the
complex conjugate. Taking the Fourier transform of \eqref{a-we}, we
obtain the Helmholtz equation:
\begin{equation}\label{a-he}
\Delta u(\textbf{r},\omega)+\kappa^2(\omega)u(\textbf{r},\omega)=
-4\pi\kappa^2(\omega)P(\textbf{r},\omega) ,
\end{equation}
where $\kappa(\omega)=\omega/c$ is the wavenumber.

The polarization may be expanded in powers of the electric field. In principle
the expansion involves infinitely many terms, but only the first few terms are
of practical importance if the nonlinearity is weak. In this paper, we consider
linear, quadratically nonlinear, and cubically nonlinear media. 

(i) A medium is linear if
\begin{equation}\label{a-lp}
P(\textbf{r},\omega)=\chi^{(1)}(\textbf{r},\omega)u(\textbf{r},\omega),
\end{equation}
where the coefficient $\chi^{(1)}(\textbf{r},\omega)$ is the first-order
susceptibility. Combining \eqref{a-he} and \eqref{a-lp} gives
\begin{equation}\label{ahe-lp}
 \Delta u(\textbf{r},
\omega)+\kappa^2(\omega)(1+4\pi\chi^{(1)}({\bf r},
\omega))u(\textbf{r},\omega)=0. 
\end{equation}

(ii) A medium is quadratically nonlinear if
\begin{equation}\label{a-qp}
P(\textbf{r},\omega)=\chi^{(1)}(\textbf{r},\omega)u(\textbf{r},\omega) +
\sum_{\omega_1+\omega_2=\omega}
\chi^{(2)}(\textbf{r},\omega_1,\omega_2) u(\textbf{r},\omega_1)u(\textbf{r},
\omega_2),
\end{equation}
where $\chi^{(2)}(\textbf{r},\omega_1,\omega_2)$ are the second-order
susceptibilities. The summation indicates that the electric fields at the 
frequencies $\omega_1$ and $\omega_2$ contribute to the
polarization at the frequency $\omega$ if $\omega_1+\omega_2=\omega$.
Second-order nonlinear effects include second harmonic generation, which is
excited by a monochromatic incident field of frequency $\omega$ in a
quadratically nonlinear medium. We assume that the nonlinear susceptibilities
are sufficiently weak, i.e., 
\[
\sum_{\omega_1+\omega_2=\omega}\chi^{(2)}(\textbf{r},\omega_1,\omega_2)u(\textbf
{r},\omega_1)u(\textbf{r},\omega_2)\ll\chi^{(1)}(\textbf{r},\omega)u(\textbf{r},
\omega).
\]
We also assume that the second-order susceptibilities have full permutation
symmetry, i.e., 
\[
\chi^{(2)}(\textbf{r},\omega_1,\omega_2)=\chi^{(2)}(\textbf{r},\omega_2,
\omega_1).
\]
Let $\omega_1=\omega$ and $\omega_2=2\omega$. Define $\kappa_j=\omega_j/c$.
Denote by $u^{(j)}$ the field corresponding to the wavenumber $\kappa_j$. It
follows from \eqref{a-he} and \eqref{a-qp} that $u^{(j)}$ satisfies
\begin{subequations}\label{ahe-qp}
 \begin{align}
&\Delta u^{(1)}({\bf r})+\kappa^2_1 (1+4\pi\chi^{(1)}({\bf r},
\omega_1))
u^{(1)}({\bf r})=-8\pi\kappa_1^2 \chi^{(2)}({\bf r}, \omega_2,
-\omega_1)\bar{u}^{(1)}({\bf r})u^{(2)}({\bf r}),\\
&\Delta u^{(2)}({\bf r})+\kappa_2^2(1+4\pi\chi^{(1)}({\bf r},
\omega_2))u^{(2)}({\bf r})=-4\pi\kappa_2^2 \chi^{(2)}({\bf r}, \omega_1,
\omega_1)\bigl(u^{(1)}({\bf r})\bigr)^2.
 \end{align}
\end{subequations}

(iii) A medium is cubically nonlinear if
\begin{equation}\label{a-cp}
P(\textbf{r},\omega)=\chi^{(1)}(\textbf{r},\omega)u(\omega)+\sum_{
\omega_1+\omega_2+\omega_3=\omega}\chi^{(3)}(\textbf{r},\omega_1,\omega_2,
\omega_3)u(\textbf{r},\omega_1)u(\textbf{r},\omega_2)u(\textbf{r},\omega_3),
\end{equation}
where $\chi^{(3)}(\textbf{r},\omega_1,\omega_2,\omega_3)$ are the third order
susceptibilities. Materials with inversion symmetry have zero second order
susceptibilities and thus fall into this category. Third-order nonlinear effects
include third-harmonic generation. We assume that
the nonlinear susceptibilities are
sufficiently weak, i.e., 
\[
\sum_{\omega_1+\omega_2+\omega_3=\omega}\chi^{(3)}(\textbf{r},\omega_1,\omega_2,
\omega_3) u(\textbf{r},\omega_1)u(\textbf{r},\omega_2)u(\textbf{r},
\omega_3)\ll\chi^{(1)} (\textbf{r},\omega)u(\textbf{r},\omega). 
\]
In addition, we assume that third-order susceptibilities have full
permutation symmetry, i.e., 
\[
\chi^{(2)}(\textbf{r},\omega_1,\omega_2,\omega_3)=\chi^{(2)}(\textbf{r},\omega_{
{\rm p} (1)},\omega_{{\rm p} (2)},\omega_{{\rm p}(3)}),
\]
where $\{{\rm p}(1), {\rm p}(2), {\rm p}(3)\}$ is a permutation of $\{1, 2,
3\}$. Suppose that a source of frequency $\omega$ is incident upon a cubic
nonlinear medium. Let $\omega_1=\omega$ and $\omega_3=3\omega$. Define
$\kappa_j=\omega_j/c$. Denote
by $u^{(j)}$ the field corresponding to the wavenumber $\kappa_j$. It
follows from \eqref{a-he} and \eqref{a-cp} that $u^{(j)}$ satisfies
\begin{subequations}\label{ahe-cp}
 \begin{align}
\Delta u^{(1)}({\bf r})+\kappa^2_1 (1+4\pi\chi^{(1)}({\bf r}, \omega_1))
u^{(1)}({\bf r})=&-12\pi\kappa_1^2 \chi^{(3)}({\bf r}, \omega_1, \omega_1, 
-\omega_1)\bar{u}^{(1)}({\bf r})\bigl(u^{(1)}({\bf r})\bigr)^2 \notag\\
& -12\pi\kappa_1^2 \chi^{(3)}({\bf r}, \omega_3, -\omega_1,
-\omega_1)u^{(3)}({\bf r})\bigl(\bar{u}^{(1)}({\bf r})\bigr)^2,\\
\Delta u^{(3)}({\bf r})+\kappa_3^2(1+4\pi\chi^{(1)}({\bf r},
\omega_3))u^{(3)}({\bf r})=&-4\pi\kappa_3^2 \chi^{(3)}({\bf r}, \omega_1,
\omega_1, \omega_1)\bigl(u^{(1)}({\bf r})\bigr)^3.
 \end{align}
\end{subequations}

\section{Point scatterer models}

Scatterers that are small compared to the wavelength can be effectively treated
as point scatterers \cite{vcl-rmp98}. The susceptibilities of a small scatterer
located at ${\bf r}_0$ can be replaced by delta functions in the following
forms:
\begin{subequations}\label{b-pt}
\begin{align}
\chi^{(1)}({\bf r}, \omega) &= \eta^{(1)}({\bf r},\omega) \delta({\bf r} - {\bf
r}_0) , \\
\chi^{(2)}({\bf r}, \omega_1,\omega_2) &= \eta^{(2)}({\bf r},\omega_1,\omega_2)
\delta({\bf r} - {\bf r}_0) , \\
\chi^{(3)}({\bf r}, \omega_1,\omega_2, \omega_3) &= \eta^{(3)}({\bf
r},\omega_1,\omega_2, \omega_3) \delta({\bf r} - {\bf r}_0),
\end{align}
\end{subequations}
where $\eta^{(j)}$ are the effective susceptibilities of the point scatterer.
The representations \eqref{b-pt} are valid outside of the interaction region,
which is of size comparable to the physical size of the scatterer. 

(i) For a collection of linear point scatterers located at ${\bf r}_k\in\mathbb
R^2, k=1, \dots, m$, the wave equation \eqref{ahe-lp} takes the form:
\begin{equation}\label{b-lp}
 \Delta u(\textbf{r},
\omega)+\kappa^2(\omega)u(\textbf{r},\omega)= - 4\pi \kappa^2(\omega)
\sum_{k=1}^m\eta^{(1)}({\bf r}_k,\omega) u(\textbf{r}_k,\omega) \delta({\bf r} -
{\bf r}_k).
\end{equation}
Identifying $u(\textbf{r},\omega)$ with $\phi(\bf r)$ and $4\pi 
\kappa^2(\omega) \eta^{(1)}({\bf r}_k,\omega) $ with $\sigma_k$, we obtain
\eqref{flhe-tf}.

(ii) For a collection of quadratically nonlinear point scatterers located at
${\bf r}_k\in\mathbb R^2, k=1, \dots, m$, the wave equation \eqref{ahe-qp} takes
the form:
\begin{subequations}\label{b-qp}
 \begin{align}
\Delta u^{(1)}({\bf r})+\kappa^2_1 u^{(1)}({\bf r})
=& - \sum_{k=1}^m \Bigl(4\pi \kappa_1^2 \eta^{(1)}({\bf r}_k,
\omega_1)u^{(1)}({\bf r}_k) 
\notag\\
&+ 8\pi\kappa_1^2 \eta^{(2)}({\bf r}_k,
\omega_2,-\omega_1)\bar{u}^{(1)}({\bf r}_k)u^{(2)}({\bf r}_k) \Bigr)\delta( {
\bf r}-{\bf r}_k),\\
\Delta u^{(2)}({\bf r})+\kappa_2^2 u^{(2)}({\bf r})
=& - \sum_{k=1}^m \Bigl(4\pi \kappa_2^2 \eta^{(1)}({\bf
r}_k,\omega_2))u^{(2)}({\bf r}_k)
\notag\\
&+ 4\pi\kappa_2^2 \eta^{(2)}({\bf r}_k, \omega_1,\omega_1)\bigl(u^{(1)}({\bf
r}_k)\bigr)^2 \Bigr) \delta( { \bf r}-{\bf r}_k).
 \end{align}
\end{subequations}
Identifying $u^{(j)}$ with $\phi^{(j)}$, $4\pi  \kappa_1^2 \eta^{(1)}({\bf
r}_k,\omega_1) $ with $\sigma^{(1)}_{k,1}$, and $8\pi\kappa_1^2 \eta^{(2)}({\bf r}_k,
\omega_2,-\omega_1)$ with $\sigma^{(2)}_{k,1}$, we obtain (\ref{flhe-qtf}a).
Identifying $4\pi  \kappa_2^2 \eta^{(1)}({\bf r}_k,\omega_2) $ with
$\sigma^{(1)}_{k,2}$, and $4\pi\kappa_2^2 \eta^{(2)}({\bf r}_k, \omega_1,\omega_1)$
with $\sigma^{(2)}_{2,k}$, we obtain (\ref{flhe-qtf}b). Note
that $\sigma^{(1)}_{k,1}$ and $\sigma^{(1)}_{k,2}$ generically take different
values, so do $\sigma^{(2)}_{k,1}$ and $\sigma^{(2)}_{k,2}$.

(iii) For a collection of cubically nonlinear point scatterers located at ${\bf
r}_k\in\mathbb R^2, k=1, \dots, m$, the wave equation \eqref{ahe-cp} takes the
form:
\begin{subequations}\label{b-cp}
 \begin{align}
\Delta u^{(1)}({\bf r})+\kappa^2_1 u^{(1)}({\bf r})
=& - \sum_{k=1}^m \Bigl(4\pi \kappa^2_1\eta^{(1)}({\bf r}_k, \omega_1))
u^{(1)}({\bf r}_k) \notag\\
&+ 12\pi\kappa_1^2 \eta^{(3)}({\bf r}_k, \omega_1, \omega_1,
-\omega_1)\bar{u}^{(1)}({\bf r}_k)\bigl(u^{(1)}({\bf r}_k)\bigr)^2 \notag\\
& +12\pi\kappa_1^2 \eta^{(3)}({\bf r}_k, \omega_3, -\omega_1,
-\omega_1)u^{(3)}({\bf r}_k)\bigl(\bar{u}^{(1)}({\bf r}_k)\bigr)^2\Bigr) \delta(
{ \bf r}-{\bf r}_k),\\
\Delta u^{(3)}({\bf r})+\kappa_3^2 u^{(3)}({\bf r})
=& - \sum_{k=1}^m \Bigl(4\pi \kappa_3^2 \eta^{(1)}({\bf r}_k,\omega_3))
u^{(3)}({\bf r}_k) \notag\\
&+ 4\pi\kappa_3^2 \eta^{(3)}({\bf r}_k, \omega_1,\omega_1,
\omega_1)\bigl(u^{(1)}({\bf r}_k)\bigr)^3\Bigr) \delta( { \bf r}-{\bf r}_k).
 \end{align}
\end{subequations}
Identifying $u^{(j)}$ with $\phi^{(j)}$, $4\pi  \kappa_1^2 \eta^{(1)}({\bf
r}_k,\omega_1) $ with $\sigma^{(1)}_{k,1}$, $12\pi\kappa_1^2 \eta^{(3)}({\bf r}_k,
\omega_1, \omega_1,-\omega_1)$ with $\sigma^{(3)}_{k,1}$, and
$12\pi\kappa_1^2 \eta^{(3)}({\bf r}_k, \omega_3,-\omega_1,- \omega_1)$ with the
second $\sigma^{(3)}_{k,2}$, we obtain (\ref{flhe-ctf}a). Identifying and $4\pi 
\kappa_3^2 \eta^{(1)}({\bf r}_k,\omega_3) $ with $\sigma^{(1)}_{k,2}$ and
$4\pi\kappa_3^2 \eta^{(3)}({\bf r}_k, \omega_1, \omega_1,\omega_1)$ with
$\sigma^{(3)}_{k,3}$ , we obtain (\ref{flhe-ctf}b).  Note that $\sigma^{(1)}_{k,1}$ and
$\sigma^{(1)}_{k,2}$ generically take different values, and so do
$\sigma^{(1)}_{k,1}$, $\sigma^{(1)}_{k,2}$ and $\sigma^{(1)}_{k,3}$.

\end{document}